\documentclass{article}
\usepackage{amsthm}
\usepackage{amsfonts}
\usepackage{amssymb}
\usepackage{amsgen}
\usepackage{amsmath}
\usepackage{amsopn}
\usepackage{verbatim}
\usepackage{xypic}
\usepackage{xspace}
\usepackage{multicol}
\usepackage{makeidx}
\usepackage{eepic}
\usepackage{upref}

\theoremstyle{plain}
\newtheorem{thm}{Theorem}[section]
\newtheorem{lem}[thm]{Lemma}
\newtheorem{prop}[thm]{Proposition}
\newtheorem{cor}[thm]{Corollary}
\newtheorem{conj}[thm]{Conjecture}
\theoremstyle{definition}
\newtheorem{prob}[thm]{Problem}

\newtheorem{rem}[thm]{Remark}
\newtheorem{rems}[thm]{Remarks}

\newcommand{\ra}{\rightarrow}
\newcommand{\ors}{{\mathbf{s}}}
\newcommand{\ort}{{\mathbf{t}}}

\newcommand{\Z}{{\mathbb Z}}

\newcommand{\N}{{\mathbb N}}
\newcommand{\ZN}{{\mathbb Z_{\ge 0}}}

\DeclareMathOperator{\pari}{{\mathfrak p}}
 \DeclareMathOperator{\sh}{{Shfl}}
\DeclareMathOperator{\re}{{Re}}

\newcommand{\gl}{{\lambda}}

\newcommand{\gk}{{\kappa}}

\newcommand{\gs}{{\sigma}}

\begin{document}

\title{Multiple Harmonic Sums I:\\
Generalizations of Wolstenholme's Theorem}
\author{Jianqiang Zhao\footnote{Email: zhaoj@eckerd.edu. This work is partially supported by
NSF grant DMS0139813 and DMS0348258}}
\date{}
\maketitle
\begin{center}
{\large Department of Mathematics, Eckerd College, St. Petersburg,
Florida 33711, USA}
\end{center}

\textbf{Abstract.} In this note we will study the $p$-divisibility
of multiple harmonic sums. In particular we provide some
generalizations of the classical Wolstenholme's Theorem to both
homogeneous and non-homogeneous sums. We make a few conjectures at
the end of the paper and provide some very convincing evidence.


\section{Introduction}
The Euler-Zagier multiple zeta functions of length $l$ are nested
generalizations of Riemann zeta function. They are defined as
\begin{equation}\label{zeta}
\zeta(\ors)=\zeta(s_1,\dots, s_l)=\sum_{0<k_ 1<\dots<k_l}
k_1^{-s_1}\cdots k_l^{-s_l}
\end{equation}
for complex variables $s_1,\dots, s_l$ satisfying
$\re(s_j)+\dots+\re(s_l)>l-j+1$ for all $j=1,\dots,l$. We call
$|\ors|=s_1+\dots+s_l$ the weight and denote the length by
$l(\ors)$. The special values of multiple zeta functions at
positive integers have significant arithmetic and algebraic
meanings, whose defining series \eqref{zeta} will be called {\em
MZV series}, including the divergent ones like $\zeta(\dots,1)$.
These obviously generalize the notion of harmonic series whose
weight is equal to 1.

MZV series are related to many aspects of number theory. One of the
most beautiful computations carried out by Euler is the following
evaluation of zeta values at even positive integers:
$$2\zeta(2m)=(-1)^{m+1}\frac{(2\pi)^{2m}}{(2m)!} B_{2m},$$
where $B_k$ are Bernoulli numbers defined by the Maclaurin series
$t/(e^t-1)= \sum_{k=0}^\infty B_k t^k/k!.$ In this paper, we will
study partial sums of MZV series which turn out to be closely
related to Bernoulli numbers too. These sums have been called
(non-alternating) {\em multiple harmonic sums} and studied by
theoretical physicists (see \cite{bl1,bl2} and their references).
Their main interest is fast computation of their exact values and
the algebraic relations among them. Our focus, however, is on
$p$-divisibility of these sums for various primes $p$ and special
attention will be paid to the cases where the sums are divisible
by higher powers of primes than ordinarily expected which is often
related to the irregular primes, i.e., primes $p$ which divide
some Bernoulli numbers $B_t$ for some positive even integers
$t<p-2$. We say in this case $(p,t)$ is an irregular pair.

\subsection{Wolstenholme's Theorem}
Put $\ors:=(s_1,\dots, s_l)\in \N^l$ and denote the $n$th partial
sum of MZV series by
\begin{equation}\label{equ:defnH}
H(\ors;n)=H(s_1,\dots,s_l;n):=\sum_{1\le k_ 1<\dots<k_l\le n}
k_1^{-s_1}\cdots k_l^{-s_l},\quad n\in \ZN.
\end{equation}
By convention we set $H(\ors;r)=0$ for $r=0,\dots,l-1$, and
$H(\emptyset;0)=1$. To facilitate our study we also define (cf.
\cite{H1,H2})
\begin{equation}\label{equ:defnS}
S(s_1,\dots,s_l;n):=\sum_{1\le k_ 1\le \dots \le k_l\le n}
k_1^{-s_1}\cdots k_l^{-s_l},\quad n\in \ZN.
\end{equation}
To save space, for an ordered set $(e_1,\dots,e_t)$ we denote by
$\{e_1,\dots, e_t\}^d$ the set formed by repeating
$(e_1,\dots,e_t)$ $d$ times and $1^d=\{1\}^d$. For example
$H(\{s\}^l;n)$ is just a partial sum of the nested zeta value
series $\zeta(s)$ of length $l$ which we refer to as a {\em
homogeneous} partial sum. The partial sums of nested harmonic
series are related to Stirling numbers $St(n,j)$ of the first kind
which are defined by the expansion
$$
x(x+1)(x+2)\cdots(x+n-1)=\sum_{j=1}^n St(n,j) x^j.
$$
We have $St(n,n)=1$, $St(n,n-1)=n(n-1)/2$, and $St(n,1)=(n-1)!$.
It's also easy to see that
\begin{equation}\label{Stirling}
St(n,j) = (n-1)! \cdot H(1^{j-1}; n-1), \text{ for} j=1,\cdots, n.
\end{equation}
Now we restrict $n$ to prime numbers.
\begin{thm}\label{Wols}
{\em (Wolstenholme 1862 \cite[p.89]{HW})} For any prime number
$p\ge 5$, $St(p,2)\equiv   0 \pmod{p^2}.$
\end{thm}
We find on the Internet the following generalization of the above theorem
by Bruck \cite{Bruck} although no proof is given there. Denote by
$\pari(m)$ the parity of $m$ which is $1$ if $m$ is odd and
$2$ if $m$ is even.
\begin{thm}\label{thm:Bruck}
For any prime number $p\ge 5$ and positive integer $l=1,\dots,
p-3,$ we have
\begin{equation*}
St(p,l+1) \equiv 0,\quad H(1^l; p-1)\equiv 0 \pmod{p^{\pari(l+1)}
}.
\end{equation*}
\end{thm}
This of course implies that $p|St(p,l)$ for $1<l<p$ which was
known to Lagrange \cite[p.87]{HW}.

It is noticed that not only $p^2$ but also $p^3$ possibly divides
$St(p,2)$, though rarely, and therefore $p^3$ possibly divides the
numerator of $H(1; p-1)$ written in the reduced form. So far we
know this happens only for $p=16843$ and $p=2124679$ among all the
primes up to 12 million (see \cite{LLV, Wagstaff, EM, BCEM,
BCEMS}). The reason, Gardiner told us in \cite{Gardiner}, is that
these two primes are the only primes $p$ in this range such that
$p$ divides the numerator of $B_{p-3}$. Bruck \cite{Bruck} further
gave a heuristic argument to show that there should be infinitely
many primes $p$ such that $p^3$ divides $S(p,2)$, which is
equivalent to say that there are infinitely many irregular pairs
$(p,p-3)$.

There are other generalizations of Wolstenholme 's Theorem to
multiple harmonic sums. Bayat proved (corrected version, see
Remark~\ref{rem:cor})
\begin{thm}\label{thm:Bayat} {\em (\cite[Thm.~3]{bayat})}
For any positive integer $s$ and prime number $p\ge k+3$ we have
\begin{equation*}
H(s; p-1)\equiv 0 \pmod{p^{\pari(s+1)} }.
\end{equation*}
\end{thm}

In \cite{sla} Slavutskii showed
\begin{thm} \label{thm:sla}
{\em (\cite[Thm.~2]{sla})} Let $s$ and $m$ be two positive
integers such that $(m,6)=1$. Let $t=(\phi(m)-1)s$ and $A_n(m)=\prod_{p|m}(1-p^{n-1})B_n$
be the Agoh's function, where $\phi(m)$ is the Euler's Phi function. Then
$$H^*(s;m-1)=\sum_{\substack{1\le k<m\\ (k,p)=1}}\frac{1}{k^s}
\equiv \begin{cases}
m A_t(m) &\pmod {m^2} \quad \text{ if $s$ is even},\\
(t/2) m^2 A_{t-1}(m)  &\pmod {m^2}  \quad \text{ if $s$ is odd}.
\end{cases}$$
\end{thm}

\subsection{Homogeneous multiple harmonic sums}
The classical result of Wolstenholme is the original motivation of
our study. We will prove the following generalization of
Thm.~\ref{Wols} and Thm.~\ref{thm:Bruck} to homogeneous multiple
harmonic sums
 (see Thm.~\ref{thm:p-1})
\begin{thm}\label{genWols}
Let $s$ and $l$ be two positive integers.  Let $p$ be an odd prime
such that $p\ge l+2$ and $p-1$ divides none of $ls$ and $ks+1$ for
$k=1,\dots, l$. Then
$$H(\{s\}^l; p-1) \equiv 0 \pmod{p^{\pari(ls-1)}}. $$
In particular, the above is always true if $p\ge ls+3$.
\end{thm}
We also look at some cases when the congruences hold modulo higher
powers of $p$. Recently, Zhou and Cai \cite{zc} prove that
\begin{thm} \label{thm:zcIn}
Let $s$ and $l$ be two positive integers. Let $p$ be a prime such
that $p\ge ls+3$
\begin{equation*} H(\{s\}^l;p-1)\equiv S(\{s\}^l;p-1)\equiv \left\{
\begin{aligned}
 (-1)^l  \frac{s(ls+1)p^2}{2(ls+2)}B_{p-ls-2}   & \pmod{p^3} \quad  \text{if}\;2\nmid ls,\\
 (-1)^{l-1} \frac{sp}{ls+1}B_{p-ls-1}\quad\  &\pmod{p^2} \quad   \text{if}\;2\,|\,ls.
\end{aligned}
\right.
\end{equation*}
\end{thm}
We will prove an analog of this in the non-homogeneous even weight
length two case (see Thm.~\ref{thm:nonhomol=2}).

\subsection{Non-homogeneous multiple harmonic sums}
The third section of this paper deals with non-homogeneous
multiple harmonic sums. We consider the length 2 case in
Thm.~\ref{thm:length2} whose proof relies heavily on generating
functions of the Bernoulli polynomials and properties of Bernoulli
numbers such as Claussen-von Staudt Theorem.
\begin{thm}
Let $s_1,s_2$ be two positive integers
and $p$ be an odd prime. Let $s_1\equiv m, s_2\equiv n\pmod {p-1}$
where $0\le m,n\le p-2$. If $m,n\ge 1$ then
\begin{equation*}
H(s_1,s_2;p-1)\equiv
\begin{cases}
\displaystyle{\frac{(-1)^n }{m+n}}{m+n\choose m}
 B_{p-m-n}  & \pmod{p}\quad \text{ if } p\ge m+n,\\
0\quad &\pmod{p}\quad \text{ if } p<m+n.
\end{cases}
\end{equation*}
\end{thm}
The same idea but more complicated computation enables us to deal
with the length 3 odd weight case completely (see
Thm.~\ref{thm:oddlength3}).
\begin{thm}
Let $p$ be an odd prime. Let $(s_1,s_2,s_3)\in \N^3$ and
$0\le l,m,n\le p-2$ such that $s_1\equiv l, s_2\equiv m,s_3\equiv n\pmod{p-1}$.
If $l,m,n\ge 1$ and $w=l+m+n$ is an {\em odd} number then
\begin{equation*}
H(l,m,n;p-1)\equiv I(l,m,n)-I(n,m,l) \pmod{p}
\end{equation*}
where $I$ is defined as follows. Let $w'=w-(p-1)$ if
$p<w<2p$ and $w'=w$ otherwise. Then
\begin{equation*}
I(l,m,n)=
\begin{cases}
0        &\text{ if }w\ge 2p,\text{ or if $l+m<p$ and $p<w<2p-1$},\\
1/2n  &\text{ if }w=p, 2p-1,\\
\displaystyle{ (-1)^{n+1}{w'\choose n} \frac{B_{p-w'}}{2w'}} & \text{otherwise}.
\end{cases}
\end{equation*}
\end{thm}

However, it seems to be extremely difficult to adopt the same
machinery for general larger length cases. For the even weight
cases in length 3, we are only able to determine the
$p$-divisibility for $H(4,3,5;p-1)$, $H(5,3,4;p-1)$ and the three
multiple harmonic sums of weight 4: $H(1,1,2;p-1)$,
$H(1,2,1;p-1)$, and $H(2,1,1;p-1)$, which are distinctly different
from the behavior of others.

Recently, M.~Hoffman studies the same kind  of questions
independently in \cite{H1,H2} from a different viewpoint. We
strongly encourage the interested reader to compare his results to
ours. For example, Hoffman defines the convolution operation on
composite of indices (see \cite[\S6]{H1}). We can apply this to a
few of the above results to find more Wolstenhomles' type
congruence in section~\ref{sec:conv}.
\begin{thm} \label{thm:convintro}
Let $p$ be a prime and $\ors\in \N^l$. Assume $p>|\ors|+2$. Then
$$H(\ors;p-1)\equiv S(\ors;p-1)\equiv 0\pmod{p}$$
provided $\ors$ has one of the following forms:
\begin{enumerate}

\item  $\ors=\big (1^m,2,1^n\big )$ for $m,n\ge 0$ and $m+n$ is
even.

\item   $\ors= \big(1^{n},2,1^{n-1},2,1^{n+1}\big )$ where $n\ge
2$ is even.

\item  $\ors= \big(1^{n+1},2,1^{n-1},2,1^{n}\big)$ where $n\ge 2$
is even.

\item   $\ors= \big(1^{n},2,1^n,2,1^n\big )$ where $n\ge 0$.
\end{enumerate}
\end{thm}

In \cite{pmod} we will look at the mod $p$ structure of the
multiple harmonic sums for lower weights.

One can also investigate the multiple harmonic sum $H(\ors; n)$
with fixed $\ors$ but varying $n$. We will carry this out in the
second part of this series \cite{2ndpart}. Such a study for
harmonic series was initiated systematically by Eswarathasan and
Levine \cite{EL} and Boyd \cite{Boyd}, independently.

The theory of Bernoulli numbers and irregular primes has a long
history, and results in this direction are scattered throughout
the mathematical literature for almost three hundred years
starting with the posthumous work ``Ars Conjectandi'' (1713) by
Jakob Bernoulli (1654-1705), see \cite{bern}. Without attempting
to be complete, we only list some of the modern references at the
end. In particular, I learned a lot from the work by Buhler,
Crandall, Ernvall, Johnson, Mets\"ankyl\"a, Sompolski,
Shokrollahi, \cite{BCEM, BCEMS, BCS, CB, EM, EM2, J}, and Wagstaff
\cite{Wagstaff} on finding irregular primes. For earlier history,
one can consult \cite{Wagstaff} and its references. Often in my
computation I use the table for irregular primes less than 12
million available online at ftp://ftp.reed.edu/users/jpb
maintained by Buhler. My interest on multiple harmonic sums was
aroused by the work of Boyd \cite{Boyd}, Bruck \cite{Bruck},
Eswarathasan and E. Levine \cite{EL}, and Gardiner \cite{Gardiner}
on the nice and surprising relations between partial sums of
harmonic series and irregular primes. I'm indebted to all of them
for their efforts on improving our knowledge of this beautiful
part of number theory.

\section{Generalizations of Wolstenholme's Theorem}
It's known to every number theorist that for every odd prime
$p$ the sum of reciprocals
of 1 to $p-1$ is congruent to 0 modulo $p$.
However, it's a little surprising (at least for me) to know that
the sum actually is congruent to 0 modulo $p^2$ if $p\ge 5$.
This remarkable theorem was proved by Wolstenholme in 1862.

\subsection{Generalization to zeta-value series}
To generalize Wolstenholme's Theorem  we need the classical
Claussen-von Staudt Theorem on Bernoulli numbers
(see, for example, \cite[p. 233, Thm.~3]{IR}):

\begin{lem}\label{Bernoullip}
For $m\in \N$, $B_{2m}+\sum_{p-1|2m} 1/p$ is an integer.
\end{lem}

We begin with a special case of our generalization which only
deals with zeta-value series, i.e., MZV series of length 1. The
general case will be built upon this.

\begin{lem}\label{lem:l=1}
Let $p$ be an odd prime and $s$ be a positive integer. Then
$$
H(s; p-1)\equiv
\begin{cases}
0 &\pmod{p^{\pari(s+1)} }\quad\mbox{if $p-1\nmid s,s+1$},\\
-1 &\pmod{\,p\,\,}\quad\mbox{if  $p-1| s$},\\
-p(n+1)/2 &\pmod{p^2} \quad\mbox{if $s+1=n(p-1)$}.
\end{cases}
$$
\end{lem}
\begin{rems}\label{rem:cor}

(1) The conditions in
Bayat's generalization of Wolstenholme's Theorem \cite[Thm.~3]{bayat}
should be corrected. For example, taking $k=2$ and $p=5$ in \cite[Thm.~3(i)]{bayat}
we only get $H(3;4)=2035/1728\not \equiv 0\pmod{5^2}$. In general, if $2k=p-1$
in \cite[Thm.~3(i)]{bayat} we find that $H(2k-1; p-1)\equiv -p \pmod{p^2}$
by taking $s=2k-1$ and $n=1$ in our lemma.

(2) Most of the lemma follows from Slavutskii's result Thm.~\ref{thm:sla}.
However, we need a more direct proof which we will reference later.
\end{rems}

\begin{proof}
If $p-1|s$ then $H(s; p-1)\equiv p-1 \equiv -1 \pmod{p}$ by
Fermat's Little Theorem. So we assume $p-1\nmid s$. This implies
that the map $a\ra a^{-s}$ is nontrivial on $(\Z/p\Z)^\times$, say
$b^s\not \equiv 1 \pmod{p} $for some $1<b<p$ . Then it's not hard
to see that $(1-b^{-s})H(s; p-1)\equiv 0 \pmod{p}$ and therefore
$H(s; p-1)\equiv 0 \pmod{p}$. This holds for any $s$ such that
$p-1\nmid s$, whether it is even or odd.

The last case of modulus $p^2$ for odd $s$ can be handled by the
same argument as in the proof of Wolstenholme's Theorem. We
produce two proofs below for both completeness and later
reference.

Let $s$ be an odd positive integer.
Choose $n$ large enough so that $t:=np(p-1)-s\ge 3$ is odd. Then
by the general form of Fermat's Little Theorem
\begin{equation*}
H(s; p-1)=\sum_{k=1}^{p-1} \frac{1}{k^{s} }
\equiv \sum_{k=1}^{p-1} \frac{k^{np(p-1)}}{k^{s} }
\equiv\sum_{k=1}^{p-1} k^t. \pmod{p^2}
\end{equation*}
By a classical result of sums of powers (see \cite[p. 229]{IR}) we
know that
\begin{equation}\label{key}
 \sum_{k=1}^{p-1} k^t =\frac{1}{t+1} (B_{t+1}(p)-B_{t+1})
     \quad\text{ for }t\ge 1,
\end{equation}
where $B_m(x)$ are the Bernoulli polynomials. Further,
$$B_{t+1}(p)=\sum_{j=0}^{t+1} {t+1 \choose j} B_j p^{t+1-j}.$$
Observing that $p B_j$ is always $p$-integral by
Lemma~\ref{Bernoullip} we have
\begin{equation}\label{pBt}
\sum_{k=1}^{p-1} k^t \equiv p B_t + \frac{t}{2} p^2 B_{t-1}
\pmod{p^2}.
\end{equation}
When $p-1\nmid s+1$ the lemma follows from the facts that
$B_j=0$ if $j>2$ is odd and that $B_{t-1}=B_{np(p-1)-s-1}$
is $p$-integral by Lemma~\ref{Bernoullip}. If $p-1|s+1$
and $s+1=m(p-1)$ then we choose $n=m$ in the above argument.
Then $t=mp(p-1)-m(p-1)+1$ is odd. So $B_t=0$, $p-1|t-1$ and
$pB_{t-1}\equiv -1\pmod{p}$ by Lemma~\ref{Bernoullip} again.
 From congruence \eqref{pBt} we get
$$H(s; p-1)\equiv -pt/2\equiv -p(n+1)/2 \pmod{p^2}.$$

In fact, there's a shorter proof for the odd case which is
not as transparent as the above proof. By binomial expansion
we see that
$$
2\sum_{k=1}^{p-1} \frac{1}{k^s}=\sum_{k=1}^{p-1}
\left(\frac{1}{k^s}+\frac{1}{(p-k)^s}\right) \equiv
\sum_{k=1}^{p-1} \frac{spk^{s-1}}{k^s(p-k)^s} \pmod{p^2}.
$$
Noticing that $1/(p-k)^s\equiv -1/k^s \pmod{p}$ we have
$$
2\sum_{k=1}^{p-1} \frac{1}{k^s}\equiv\sum_{k=1}^{p-1}
\frac{-sp}{k^{s+1}}\equiv 0 \pmod{p^2}
$$
whenever $p-1\nmid s+1$ because $p$ divides
$\sum_{k=1}^{p-1} 1/k^{s+1}$ by the even case.
\end{proof}

When $s=p^e$ we can work more carefully with binomial expansion in the
shorter proof and see that $p^{2+e}$ divides $H(s; p-1)$. When
$e=1$ this explains the fact that 125 divides $H(5;4)$. We record
the phenomenon in the following proposition.

\begin{prop}\label{prop:Hp-1}
Let $s$ be a positive integer and $v_p(s)=v$ and $v_p(s+1)=u$
(so that either $u$ or $v$ is 0).
If $p\ge 5$ is a regular prime then
\begin{equation*}
v_p\big(H(s;p-1)\big)=
\begin{cases}
0 &\text{if }p-1|s,\\
v+1 &\text{if $s$ is even and $p-1\nmid s$},\\
u+v+1 &\text{if $s$ is odd and $p-1| s+1$},\\
u+v+2 &\text{if $s$ is odd and $p-1\nmid s+1$}.
\end{cases}
\end{equation*}
Suppose $p$ is irregular and let $1\le m<p(p-1)$ such that
$$m\equiv \begin{cases}
-s  & \pmod{p(p-1)} \text{ if $s$ is even},\\
-s-1 & \pmod{p(p-1)} \text{ if $s$ is odd}.
\end{cases}$$
If $p-1\nmid s+2,s+3$ and $p^2\nmid B_m/m$ then
the nonzero valuations can increase by at most $1$.
\end{prop}
\begin{proof} We only consider the case when $p-1\nmid s$.
Suppose $s=p^va$ and $p\nmid a$.
Let $e\ge v+2$ be any positive integer
such that $t:=p^e(p-1)- s>1$. Then $p-1\nmid t$.
It follows from Fermat's Little Theorem that
$$H(s;p-1)\equiv  \sum_{k=1}^{p-1} k^t  \pmod{p^{v+3}}.$$
If $s$ is even then $t$ is even and
\begin{alignat*}{3}
H(s;p-1)\equiv & pB_t+ \frac{t(t-1)}{6}p^3B_{t-2} &  \pmod{p^{v+3}}\\
\equiv & -a p^{v+1} \frac{B_{t'}}{t'}
+ \frac{(t-1)(t-2)}{6 t''}p^{v+3} B_{t''}  & \pmod{p^{v+3}}
\end{alignat*}
by Kummer congruences, where $t\equiv t', t-2\equiv t'' \pmod{p(p-1)}$
and $1\le t',t''< p(p-1)$. Note that
$v_p$-valuation of the first term is $v+1$ or higher depending
on whether $(p,t')$ is regular pair or not.
The smallest $v_p$-valuation
of the second term is $v+2$ which happens if and only if
$p-1|t''$. If $p$ is irregular then $p-1\nmid t''$ by our
assumption $p-1\nmid s+2$ and hence the second term is always
divisible by $p^{v+3}$. It follows that the $v_p$-valuation
is $v+1$ if $p$ is regular and is at most $v+2$ if $p$ is
irregular since we assumed $p^2\nmid B_{t'}/t'.$ This proves the
proposition when $s$ is even.

When $s$ is odd, then we need to consider two cases: $u=0$
or $v=0$. Both proofs in these two cases are similar to the
even case and hence we leave the details to the interested
readers.
\end{proof}
\begin{rem} \label{3cons}
One can improve the above by a case by case analysis modulo higher
$p$-powers. For example, one should be able to prove that the
nonzero valuations can increase by at most 6 if $p$ is irregular
less than 12 million.
\end{rem}

\begin{rem} \label{rem:highcon}
Numerical evidence shows that if $p$ is irregular then the
nonzero valuations seems to increase by at most 1 in most cases.
The first counterexample appears with $p=37$, and $s=1048$. Note
that $s$ is even and $p-1\nmid s$ so that if one leaves out the
conditions in the lemma then the prediction would say
$v_p(H(s;p-1))$ is at most 2 because 37 is irregular. But we have
$v_{37}(H(1048;36))=3$ because $37^2|B_{p(p-1)-s}=B_{284}$.
\end{rem}

\begin{cor} \label{cor:sbig7}
Let $s\ge 4$ be a positive integer. Let $p\ge 3$ be a prime.
If $p$ is irregular then we assume it satisfies the conditions
in the proceeding proposition. Then $H(s;p-1)\not \equiv 0 \pmod{p^s}$.
\end{cor}
\begin{proof} The case $p=3$ can be proved directly because $1+2^s<3^s$
if $s\ge 2$ and
$$1+\frac{1}{2^s}=\frac{1+2^s}{2^s}\not\equiv 0\pmod {3^s}.$$
Suppose $p\ge 5$. Let $s=p^v a$ where $p\nmid a$. If $p$ is regular
then by Prop.~\ref{prop:Hp-1} the largest value of $v_p(H(s;p-1)$ is
$v+2$ which is less than $s$ because $s\ge 4$ and $p\ge 5$.
If $p$ is irregular satisfying the conditions in Prop.~\ref{prop:Hp-1}
then the largest value of $v_p(H(s;p-1)$ is at most $v+3$ which is still
less than $s$ because $s\ge 4$ and $p\ge 37$.
\end{proof}

\subsection{$p$-divisibility, Bernoulli numbers, and irregular primes}
Numerical evidence shows that congruences in Lemma~\ref{lem:l=1}
is not always optimal. For any zeta value series, every once in a
while, a higher than expected power of $p$ divides its $(p-1)$-st
partial sum. A closer look of this phenomenon reveals that all
such primes are irregular primes. Going through the proof of
Lemma~\ref{lem:l=1} a bit more carefully one can obtain the
following improvement.

\begin{thm}\label{thm:bernoulli}
Suppose $n$ is a positive integer and $p$ is an odd prime
such that $p\ge 2n+3$. Then we have the congruences:
\begin{equation}\label{allcong}
\frac{-2}{2n-1}\cdot H(2n-1;p-1)\equiv p\cdot H(2n;p-1) \equiv
p^2\cdot \frac{2n}{2n+1}\cdot B_{p-2n-1}
\pmod{p^3}.
\end{equation}
Therefore the following statements are equivalent:
\begin{alignat*}{3}
(1). & \quad B_{p-2n-1}\equiv 0&\pmod{\, p\, }.\\
(2). & \quad  H(2n;p-1)=\sum_{k=1}^{p-1} 1/k^{2n} \equiv 0&\pmod{p^2}.\\
(3). & \quad H(2n-1;p-1)=\sum_{k=1}^{p-1}  1/k^{2n-1} \equiv 0&\pmod{p^3}.\\
(4). & \quad H(n,n;p-1)=\sum_{1\le k_1<k_2<p} 1/(k_1k_2)^n
    \equiv 0 &\pmod{p^2}.
\end{alignat*}
\end{thm}

\begin{proof}
The congruence relation \eqref{allcong} and the equivalence of (1)
to (3) in the theorem follows from the two congruences in
\cite{ELeh} after \cite[(16)]{EL}. See also \cite[p.~281]{Gl}. The
equivalence of (2) and (4) follows immediately from the shuffle
relation\footnote{Some authors call this a stuffle relation, i.e.,
``shuffling'' plus ``stuffing''.}:
\begin{equation}\label{nodduse}
H(n;p-1)^2= 2 H(n,n;p-1)+ H(2n;p-1).
\end{equation}
\end{proof}
\begin{rem} This theorem also follows directly from a result of
\cite{zc}. See Thm.~\ref{thm:zc}.
\end{rem}

Gardiner \cite{Gardiner} proves the special case of the equivalence
when $n=1$. He has one more equivalence condition which involves the
combinatorial number $2p \choose p$. Other variations of the
classical Wolstenholme's Theorem can be found in \cite{ALkan}.
It would be an interesting problem to find the analog for
the zeta value series. It is also worth
mentioning that it is not known whether
${2n \choose n} \equiv 2 \pmod{n^3}$ would imply
$n$ is a prime, which is called the converse of
Wolstenholm's Theorem \cite{mc}.

\subsection{Hoffman's convolution}\label{sec:Hoff}
In this section we recall Hoffman's convolution operation of the
composites and provide new congruence modulo a prime square. These
results will be extremely useful when we deal with multiple
harmonic sums of arbitrary lengths.

Let's first recall the definitions. Let $k$ be a positive integer
and $\ors=(i_1,\dots,i_k)$ of weight $n=|\ors|$. We define the
power set to be the partial sum sequence of $\ors$:
$P(\ors)=(i_1,i_1+i_2,\dots, i_1+\cdots+i_{k-1})$ as a subset of
$(1,2,\dots,n)$. Clearly $P$ provides a one-to-one correspondence
between the composites of weight $n$ and the subsets of
$(1,2,\dots,n-1)$. Then $\ors^\ast$ is the composite of weight $n$
corresponding to the complimentary subset of $P(\ors)$ in
$(1,2,\dots,n-1)$. Namely,
$$\ors^\ast=P^{-1}\big((1,2,\dots,n-1)-P(\ors)\big).$$
It's easy to see that $\ors^{\ast \ast}=\ors$ so $\ast$ is indeed
a convolution. For example, if $i_1,i_k\ge 1, i_2,\dots,i_{k-1}\ge
2$ then we have
$$(i_1,\dots,i_k)^\ast=(1^{i_1-1},2,1^{i_2-2},2,1^{i_3-2},
 \dots,2,1^{i_{k-1}-2},2,1^{i_k-1}).$$
Further we set the reversal $\overline{\ors}=(s_l,\dots,s_1).$ The
following important result is due to Hoffman:

\begin{thm} \label{thm:Hoff} {\em (\cite[Thm.~6.8]{H1})}
For $\ors=(s_1,\dots,s_l)$ and any positive integer $n$ define
\begin{equation}\label{equ:S}
S(\ors;n)=\sum_{1\le k_1\le \dots\le k_l\le n}
 k_1^{-s_1}\cdots k_l^{-s_l}.
\end{equation}
Then for all prime $p$
\begin{alignat}{2}
 \label{equ:conv}
S(\ors^\ast;p-1)\equiv &\, -S(\ors;p-1)  &\pmod{p}, \\
S(\overline{\ors};p-1)\equiv &\, (-1)^{|\ors|} S(\ors ;p-1)
  &\pmod{p}. \label{equ:rev}
\end{alignat}
We also have the following equalities:
\begin{equation}\label{equ:SH}
S(\ors;n)=\sum_{\mathbf{i}\preceq \ors} H(\mathbf{i};n),
\end{equation}
where $\mathbf{i}\prec \ors$ means $\mathbf{i}$ can be obtained
from $\ors$ by combining some of its parts, and
\begin{equation}\label{equ:partS}
(-1)^{l(\ors)}S(\overline{\ors};p-1)= \sum_{\bigsqcup_{j=1}^l
\ors_i=\ors} (-1)^l\prod_{j=1}^l H(\ors_j;p-1) \pmod{p}
\end{equation}
where $\bigsqcup_{j=1}^l \ors_i$ is the catenation of $\ors_1$ to
$\ors_l$.
\end{thm}
\begin{proof} Note that in \cite{H1} the right hand side of \eqref{equ:S}
is denoted by $S(\overline{\ors};n)$ instead. But it's not
difficult to verify that for any $\ors$
\begin{equation}
\overline{\ors}^\ast=\overline{\ors^\ast}.
\end{equation}
So \eqref{equ:conv} is equivalent to \cite[Thm.~6.8]{H1}.
Congruence \eqref{equ:rev} follows readily from the substitution
of indices: $k\to p-k$.

Equations \eqref{equ:partS} and \eqref{equ:SH} relate $S$-version
multiple harmonic series and our $H$-version (denoted by $A$ by
Hoffman). Equation \eqref{equ:SH} is \cite[Eq.~(7)]{H1}) and
\eqref{equ:partS} is equivalent to the second unlabelled formula
in the proof of \cite[Thm.~6.8]{H1}.
\end{proof}

Next we want to provide a more precise version of congruence
\eqref{equ:conv}.
\begin{thm}\label{cong:modp2}
Let $\ors$ be any composite of weight $w$. Let $p$ be an arbitrary
odd prime. Then
\begin{equation}
\label{equ:conv1}
  -S(\ors^\ast;p-1)\equiv S(\ors;p-1)+p\left(\sum_{\ort\preceq \ors}
  H(\ort\sqcup\{1\}; p-1)\right) \pmod{p^2}.
\end{equation}
\end{thm}
\begin{proof} For any sequence $\{f(n)\}_{n\ge -1}, f(-1)=0$
we know there are two operators $\Sigma$ and $\nabla$:
$$ \Sigma f(n)=\sum_{i=0}^n f(i),\quad
 \nabla f(n)=f(n)-f(n-1), \quad \forall n\ge 0.$$
Recall that $\Sigma\nabla S(\ors;n)=-S(\ors^\ast;n)$ for any
$\ors$ and positive number $n$ by \cite[Thm.~4.2]{H2}. Then for
any odd prime $p$ we have
$$
\aligned -S(\ors^\ast;p-1)=&\Sigma\nabla
S(\ors;p-1)(n)=\sum_{i=0}^{p-1} {p \choose i+1} (-1)^i S(\ors;i) \\
=& S(\ors;p-1)(n)+p\left( \sum_{i=0}^{p-2} \frac{1}{i+1}{p-1 \choose i}(-1)^i S(\ors;i)\right) \\
\equiv & S(\ors;p-1)(n)+p\left( \sum_{i=0}^{p-2}\frac{1}{i+1} S(\ors;i)\right) &\pmod{p^2}\\
\equiv & S(\ors;p-1)(n)+p \left(\sum_{i=0}^{p-2}
 \sum_{\ort\preceq\ors}\frac{H(\ort; i)}{i+1}  \right)  &\pmod{p^2}\\
\equiv & S(\ors;p-1)(n)+p\left(\sum_{\ort\preceq \ors}
  H(\ort\sqcup\{1\}; p-1)\right)  &\pmod{p^2}.
\endaligned
$$
\end{proof}

\subsection{Wolstenholme type theorem for homogeneous sums}
In the above we have studied the $p$-divisibility of $H(s;p-1)$
for positive integers $s$. It is difficult to start with arbitrary
$\ors$ so we only consider the homogeneous multiple harmonic sums
at the moment.

First we can easily verify the following shuffle relations: for
any positive integers $n,m,s_1,\dots,s_l$,
$$
H(m;n)\cdot H(s_1,\dots,s_l;n)=\sum_{\ors\in \sh(m,
(s_1,\dots,s_l))} H(\ors;n)+ \sum_{j=1}^l
H\big(s_1,\dots,s_{j-1},s_j+m,s_{j+1},\dots,s_l;n\big)
$$
where for any two ordered sets $(r_1,\dots,r_t)$ and
$(r_{t+1},\dots,r_n)$ the shuffle operation is defined by
$$\sh\big((r_1,\dots,r_t), (r_{t+1},\dots,r_n)\big):=
\bigcup_{\substack{\gs\text{ permutes } \{1,\dots,n\}, \\
\gs^{-1}(1)<\cdots<\gs^{-1}(t),\\
\gs^{-1}(t+1)<\cdots<\gs^{-1}(n)}}
 \big(r_{\gs(1)},\dots,r_{\gs(n)}\big).$$
Hence for any $k=1,\dots, l-1,$ we have
$$H(ks;n)\cdot H\big(\{s\}^{l-k};n\big)=
\sum_{\ors\in \sh\big(\{ks\},\{s\}^{l-k}\big)} H(\ors;n )+
\sum_{\ors\in\sh\big(\{(k+1)s\},\{s\}^{l-k-1}\big)}H(\ors;n).
$$
Applying $\sum_{k=1}^{l-1} (-1)^{k-1}$ on both sides we get
\begin{lem}\label{lem:homorec}
Let $s,l$ and $n$ be three positive integers. Then
\begin{equation}\label{induction}
l H\big(\{s\}^l;n\big)=\sum_{k=1}^l (-1)^{k-1} H(ks;n) \cdot
H\big(\{s\}^{l-k};n\big).
\end{equation}
\end{lem}

We now need a formula of the homogeneous multiple harmonic sum
$H\big(\{s\}^l;n\big)$ in terms of partial sums of ordinary zeta
value series. Let $P(l)$ be the set of unordered partitions of
$l$. For example, $P(2)=\{(1,1),(2)\}, P(3)=\{(1,1,1),
(1,2),(3)\}$ and so on.

\begin{lem}\label{lem:cgl}
Let $s,l$ and $n$ be positive integers. For
$\gl=(\gl_1,\dots,\gl_r)\in P(l)$ we put $H_\gl(s;n)=\prod_{i=1}^r
H(\gl_i s;n)$ and define $c_\gl=(-1)^{r-l} \prod_{i=1}^r
(\gl_i-1)!.$ Then
\begin{equation}\label{cgl}
l! H\big(\{s\}^l;n\big)= \sum_{\gl\in P(l)} c_\gl H_\gl(s;n).
\end{equation}
\end{lem}
\begin{proof} It follows easily from \cite[Thm.~2.3]{H2}.
\end{proof}

It is obvious that $c_{(1,\cdots,1)}=1.$ When $l=3$ we have
$c_{(1,2)}= -2c_1+c_2=-3,$ which implies
\begin{equation} \label{Hsss}
6H\big(\{s\}^3;n\big) = H(s;n)^3-3H(s;n)H(2s;n)+2H(3s;n).
\end{equation}
When $l=4$ we have $c_{(1,1,2)}=c_{(1,2})-3c_{(1,1)}=-6,$
$c_{(1,3)}= c_3+6c_1=8,$ and $c_{(2,2)}=3c_2=-3,$ which implies
\begin{equation*}
24H\big(\{s\}^4;n\big) = H(s;n)^4-6H(s;n)^2H(2s;n)
+8H(s;n)H(3s;n)-3H(2s;n)^2-6H(4s;n).
\end{equation*}

\begin{thm} \label{thm:p-1}
Let $s$ and $l$ be two positive integers.  Let $p$ be an odd prime
such that $p\ge l+2$ and $p-1$ divides none of $sl$ and $ks+1$ for
$k=1,\dots,l$. Then the homogeneous multiple harmonic sum
$$
H(\{s\}^l; p-1)\equiv S(\{s\}^l; p-1)\equiv  0
\pmod{p^{\pari(ls-1)}}.
$$
In particular, if $p\ge ls+3$ then the above is always true and so
$p| H(\{s\}^l; p-1)$.
\end{thm}
\begin{proof}
Congruence for $H$ follows from Eq.~\eqref{cgl} and
Lemma~\ref{lem:l=1}. Congruence for $S$ then follows from
\eqref{equ:partS}.
\end{proof}

Recently, Zhou and Cai obtain an improved version of the above
result, see \cite{zc}
\begin{thm} \label{thm:zc}
Let $s$ and $l$ be two positive integers. Let $p$ be a prime such
that $p\ge ls+3$
\begin{equation*} H(\{s\}^l;p-1)\equiv S(\{s\}^l;p-1)\equiv
\left\{
\begin{aligned}
 (-1)^l  \frac{s(ls+1)p^2}{2(ls+2)}B_{p-ls-2}   & \pmod{p^3} \quad  \text{if}\;2\nmid ls,\\
 (-1)^{l-1} \frac{sp}{ls+1}B_{p-ls-1}\quad\  &\pmod{p^2} \quad   \text{if}\;2\,|\,ls.
\end{aligned}
\right.
\end{equation*}
\end{thm}
\begin{proof} Congruence for $H$ follows from \cite{zc}.
Congruence for $S$ then follows from \eqref{equ:partS} and an
induction on $l$.
\end{proof}

\subsection{Higher divisibility and distribution of irregular pairs}
 From Thm.~\ref{thm:zc} we immediately have
\begin{prop}\label{13rd}
Let $p$ be an odd prime. Suppose $(p, p-ls-2)$ is an irregular
pair (so both $l$ and $s$ must be odd numbers). Then
\begin{equation*}
H(\{s\}^l;p-1)\equiv S(\{s\}^l;p-1)\equiv 0 \pmod{p^3}.
\end{equation*}
Suppose $(p, p-ls-1)$ is an irregular pair. Then
\begin{equation*}
H(\{s\}^l;p-1)\equiv S(\{s\}^l;p-1)\equiv 0 \pmod{p^2}.
\end{equation*}
\end{prop}

One of the first instances of this proposition is when $s=1$ and
$l=3$ given by the first irregular pair $(37,32)$, which, by our
proposition, implies that $H(1,1,1;36)\equiv 0\pmod{37^3}.$ The
case when $s=3$ and $l=3$ appears first with the irregular pair
$(9311,9300)$. So we know $H(3,3,3;9310)\equiv 0\pmod{9311^3}$.

We believe these are just the first two of infinitely many such pairs
because evidently not only Bernoulli numbers but also the difference
$p-t$ for irregular pairs  $(p,t)$ are evenly distributed modulo
any prime. Precisely, we have the following
\begin{conj}\label{conj:dist}
For any fixed positive integer $M$ and integer $c$
such that $0\le c<M$ we have
\begin{equation}\label{density}
\lim_{X\to \infty} \frac{\sharp\{(p,t): \text{prime } p|B_t,
t \text{ even}, p<X, p-t\equiv c \pmod{M}\}}
{\sharp\{(p,t): \text{prime } p|B_t,\,
t \text { even}, p<X\}}=\begin{cases}
0  &\text{ if $2|c$, $2|M$,} \\
1/M  &\text{ if $2\nmid M$,}\\
2/M  &\text{ if $2\nmid c$, $2|M$.}
\end{cases}
\end{equation}
Further we can replace the sets by restricting $p$ to all
irregular primes with a fixed irregular index which is defined
as the number of such pairs for a fixed $p$.
\end{conj}
As a result we expect there are about one-third irregular pairs
satisfying the conditions in Prop.~\ref{13rd}. In
Table~\ref{Ta:Irr} we count the first 11,000 irregular pairs. We
denote by $N(k,m)$ the number of irregular pairs $(p,t)$
satisfying  $p-t\equiv k \pmod 3$ in the top $m$ irregular pairs,
$0\le k\le 2$, and by $P(k,m)\%$ the percentage of such pairs. For
irregular primes with fixed index we compiled some more tables in
the Appendix at the end of this paper.

\begin{table}[h]
\begin{center}
\begin{tabular}{  ||c|c|c|c|c|c|c|c|c|c|c|c|| } \hline
$m$& 1,000&2,000&3,000&4,000&5,000&6,000&7,000&8,000&9,000&10,000&11,000\\ \hline
$N(0,m)$&322&664&996&1318&1637&1978&2310&2628&2968&3310&3672  \\ \hline
$P(0,m)$&32.20&33.20&33.20&32.95&32.74&32.97&33.00&32.85&32.98&33.10&33.38 \\ \hline
$N(1,m)$&343&666&1005&1352&1676&1999&2340&2683&3017&3354&3673  \\ \hline
$P(1,m)$&34.30&33.30&33.50&33.80&33.52 &33.32&33.43&33.54& 33.52 &33.54&33.39 \\ \hline
$N(2,m)$&335&670&999&1330&1687&2023&2350&2689&3015&3336&3655\\ \hline
$P(2,m)$&33.50&33.50&33.30&33.25&33.74&33.72&33.57 &33.61& 33.50&33.36&33.22 \\ \hline
\end{tabular}
\caption{Distribution of $p-t\pmod 3$ for irregular pairs $(p,t)$.}
\label{Ta:Irr}
\end{center}
\end{table}

\section{Non-homogeneous sums}
Having dealt with homogeneous multiple harmonic sums we would like
to do some initial experiments on the non-homogeneous ones.

\subsection{Non-homogeneous sums of length 2}
We begin with $\zeta(1,2)$ and $\zeta(2,1)$ series. For any
positive integer $n$ from the shuffle relation we have
\begin{equation}\label{H12shfl}
H(1,2;n)+H(2,1;n)=H(1;n)\cdot H(2;n)-H(3;n).
\end{equation}
However, it seems to be quite difficult to disentangle $H(1,2;n)$
from $H(2,1;n)$. Maple Computation for primes up to 20000 confirms
the following
\begin{thm} \label{thm:length2}
Let $s_1,s_2$ be two positive integers
and $p$ be an odd prime. Let $s_1\equiv m, s_2\equiv n\pmod {p-1}$
where $0\le m,n\le p-2$. Then
\begin{equation*}
H(s_1,s_2;p-1)\equiv
\begin{cases}
1\quad &\pmod{p}\quad\text{ if }(m,n)=(0,0), (1,0),\\
-1\quad &\pmod{p}\quad\text{ if }(m,n)=(0,1),\\
\displaystyle{\frac{(-1)^n }{m+n}}{m+n\choose m}
    B_{p-m-n}  & \pmod{p}\quad \text{ if } p\ge m+n \text{ and  }m,n\ge 1,\\
0\quad &\pmod{p}\quad \text{ otherwise}.
\end{cases}
\end{equation*}
In particular, if $p\ge s+t$ for two positive integers $s$ and $t$
then we always have
\begin{equation*}
 S(s,t;p-1)\equiv H(s,t;p-1)\equiv \frac{(-1)^t }{s+t}
 {s+t\choose s} B_{p-s-t}  \pmod{p}.
\end{equation*}

\end{thm}
\begin{proof} We leave the trivial cases to the interested readers and
assume in the rest of the proof that $m,n\ge 1$.

Let $M=p-1-m$ and $N=p-1-n$.  Then $1\le M,N\le p-2$. By Fermat's Little Theorem
$$H(m,n;p-1)\equiv \sum_{1\le k_1<k_2\le p-1} k_1^Mk_2^N \pmod{p}.$$
Define the formal power series in two variables
$$f(x,y)=\sum_{r,s=0}^{\infty}
\left(\sum_{1\le k_1<k_2\le p-1} k_1^rk_2^s\right) \frac{x^r
y^s}{r!s!}.$$
Exchanging summation we get
\begin{align*}
f(x,y)=&\sum_{1\le k_1<k_2\le p-1} e^{k_1x+k_2y}
=\frac{e^{p(x+y)}-e^{x+y}}{(e^{x+y}-1)(e^x-1)}
    -\frac{(e^{py}-e^y)e^x}{(e^y-1)(e^x-1)} \\
=&\sum_{i=-1}^{\infty}
\sum_{j=-1}^{\infty} \frac{B_{j+1}(p)-B_{j+1}(1)}{(i+1)!(j+1)!}
 \Bigl(B_{i+1} x^i (x+y)^j-B_{i+1}(1) x^i y^j\Bigr)
\end{align*}
where $B_m(1)=B_m$  if $m\ne 1$ and $B_1(1)=-B_1=\frac 12$. We only care
about the above sum when $1\le i,j\le p-2$ since
$1\le M,N\le p-2$. So we may as well replace $B_m(1)$ by $B_m$
everywhere in the last display and throw away the terms with $j=0,\pm 1$.
The resulting power series is
\begin{equation*}
g(x,y)\equiv  \sum_{i=-1}^{\infty}   \sum_{j=2}^{\infty}
\sum_{l=1}^{j} \frac{pB_{i+1}B_j}{(i+1)!j!} {j\choose l} x^{l+i}y^{j-l}
\pmod{p}.
\end{equation*}
Let $l+i=M$ and $j-l=N.$ Then we see that the coefficient of  $x^My^N$ is
\begin{equation}\label{BB}
 \sum_{l=1}^{M}
\frac{pB_{M-l+1}B_{N+l}}{(M-l+1)!(N+l)!}{N+l\choose l}
=\sum_{l=1}^{M} \frac{pB_{M-l+1}B_{N+l}}{(M-l+1)!N!l!}.
\end{equation}
Note that $0\le M-l+1\le M\le p-2$ and $2\le N+1\le N+l\le M+N\le 2p-4$.
If $M+N<p-2$ (i.e. $m+n>p$) then \eqref{BB} is always congruent to 0 mod $p$.
Otherwise, if $M+N\ge p-2$ then all the terms are in $p\Z_p$ except
when $N+l=p-1$. Hence
$$\frac{1}{M!N!}\sum_{1\le k_1<k_2\le p-1} k_1^M k_2^N
\equiv \frac{pB_{M+N+2-p}B_{p-1}}{(M+N+2-p)!N!(p-1-N)!} \pmod{p}.$$
So finally we arrive at
$$H(m,n;p-1)=\sum_{1\le k_1<k_2\le p-1} \frac{1}{k_1^mk_2^n}
\equiv \frac{-(p-m-1)!B_{p-m-n}}{n!(p-m-n)!}  \pmod{p}.$$ One can
now use Wilson's Theorem to get the final congruence for $H$ in
our theorem without too much difficulty. For $S$ we now use
$S(s,t)=H(s,t)+H(s+t)$ and Wolstenholme's theorem.
\end{proof}

Taking $m=1,n=2$ in the theorem we obtain $H(1,2;p-1)\equiv
B_{p-3}\pmod{p}$ for $p\ge 3$. We verified this on Maple for the
only two known irregular pairs of the form $(p,p-3)$, namely,
$p=16843$ and $p=2124679$. If we take $m=2,n=3$ we find that
$H(2,3;p-1)\equiv B_{p-5}\pmod{p}$ if $p\ge 5$. There is only one
irregular pair of the form $(p,p-5)$ among all primes less than 12
million, namely, $(37,32)$. Indeed, for $p=37$ computation shows
that $H(2,3;36)\equiv 0 \pmod{37}$. As a matter of fact, we
formulated Thm.~\ref{thm:length2} only after we had found these
intriguing examples.

We now can provide an analog of Thm.~\ref{thm:zc} in the
non-homogeneous case of length 2.
\begin{thm} \label{thm:nonhomol=2}
Let $p$ be an odd prime. Suppose $s$ and $t$ are two positive
integers of same parity such that $p>s+t+1$. Then
$$\aligned H(s, t; p-1) \equiv &\,  p
\left[(-1)^s t {s+t+1\choose s}-(-1)^s s{s+t+1\choose t}
 -s-t\right] \frac{B_{p-s-t-1}}{2(s+t+1)} \pmod {p^2}\\
 S(s, t; p-1) \equiv &\, p
\left[(-1)^s t {s+t+1\choose s}-(-1)^s s{s+t+1\choose t}
 +s+t\right] \frac{B_{p-s-t-1}}{2(s+t+1)} \pmod {p^2}.
\endaligned $$
\end{thm}
\begin{proof}
By the shuffle relation (dropping $p-1$ again) we see that
\begin{equation} \label{shflHs1s2}
H(s)\cdot H(t)= H(s,t) +  H(t,s) +  H(s+t).
\end{equation}
By the conditions on $s_1$ and $s_2$ we know from \eqref{allcong}
\begin{equation} \label{s1s2cong1}
H(s)\cdot H(t)\equiv 0,  H(s+t) \equiv \frac{p(s+t)}{s+t+1}
B_{p-s-t-1} \pmod {p^2}.
\end{equation}
Therefore
\begin{equation}\label{equ:s1s2}
H(s, t )+H(t,s )\equiv \frac{-p(s+t)}{s+t+1} B_{p-s-t-1}\pmod
{p^2}.
\end{equation}
Moreover,  by the old substitution trick $i,j\to p-i,p-j$
$$
\aligned H(s,t) =& \sum_{1\le j<i<p} \frac{1}{(p-i)^s (p-j)^t}\\
 \equiv & \sum_{1\le j<i<p}  \frac{1}{i^s j^t}
     \Big(1+\frac{p}{i}\Big)^s \Big(1+\frac{p}{j}\Big)^t
 \pmod{p^2}\qquad(s+t \text{ is even})\\
 \equiv & \sum_{1\le j<i<p}  \frac{1}{i^s j^t}  \Big(1+\frac{ps}{i}+\frac{pt}{j}\Big)
 \pmod{p^2} \\
 \equiv & H(t,s)+psH(t,s+1)+ptH(t+1,s)
 \pmod{p^2} \\
 \equiv & H(t,s)+p
\left[ (-1)^{s+1}s{s+t+1\choose t}+(-1)^s t{s+t+1\choose
t+1}\right] \frac{B_{p-s-t-1}}{s+t+1}
 \pmod{p^2}.
\endaligned$$
Combined with \eqref{equ:s1s2} this completes the proof of the
congruence for $H$. Then the $S$ part follows from the identity
$S(s,t)=H(s,t)+H(s+t).$
\end{proof}

\subsection{The case of palindrome $\ors$} Recall that for
$\ors=(s_1,\dots,s_l)$ we have set its reversal
$\overline{\ors}=(s_l,\dots,s_1).$
\begin{lem} \label{lem:revorder}
Let $p$ be an odd prime. Let $l$ be a positive integer and
$\ors\in\N^l$. Let $|\ors|=\sum_{i=1}^l s_i$ be the weight of
$\ors$. Then
\begin{align*}
 H(\ors;p-1)\equiv& \,(-1)^{|\ors|}H(\overline{\ors};p-1)  \pmod{p},\\
S(\ors;p-1)\equiv & \, (-1)^{|\ors|}S(\overline{\ors};p-1)  \pmod
{p}.
\end{align*}
\end{lem}
\begin{proof} Use the old substitution trick $k_i\to p-k_i$ for all $i$
in the definitions \eqref{equ:defnH} and \eqref{equ:defnS}.
\end{proof}
An immediate consequence of this lemma is
\begin{cor}\label{cor:palin}
Let $p$ be an odd prime. If $\ors=\overline{\ors}$ and $|\ors|$ is
odd then
$$H(\ors;p-1)\equiv S(\ors;p-1) \equiv0 \pmod {p}.$$
\end{cor}
On the contrary, a lot of examples show that if the weight
$|\ors|\ge 6$ is even and if the length is bigger than 2, then we
often have $H(\ors;p-1)\not\equiv 0 \pmod {p}$ when $p$ is large
(say, $p\ge 2|\ors|$), even in the case that $\ors$ is a
palindrome. For example, in length 3 if $\ors\ne(4,3,5),(5,3,4)$
then this seems to be always the case (see
Problem~\ref{prob:H435}). A remarkable different pattern occurs
for length 3 weight 4 case which we will consider in
subsection~\ref{sec:d3wt4}. Clarifying this completely might be a
crucial step to understand the structure of $H(\ors;p-1)$ in
general.

\subsection{Multiple harmonic sums of length 3 with odd weight}
One may want to generalize Thm.~\ref{thm:length2} to multiple
harmonic sums of longer lengths. However, the proofs become much
more involved. Extensive computation confirms  the following
\begin{thm}  \label{thm:oddlength3}
Let $p$ be an odd prime. Let $(s_1,s_2,s_3)\in \N^3$ and
$0\le l,m,n\le p-2$ such that $s_1\equiv l, s_2\equiv m,s_3\equiv n\pmod{p-1}$.
Then
\begin{align*}
H(0,m,n;p-1)&\equiv H(m-1,n;p-1)-H(m,n;p-1),  &\pmod{p}, \\
H(l,0,n;p-1)&\equiv H(l,n-1;p-1)-H(l,n;p-1)-H(l-1,n;p-1)  &\pmod{p}, \\
H(l,m,0;p-1)&\equiv -H(l,m-1;p-1)-H(l,m;p-1)&\pmod{p} .
\end{align*}
Suppose $l,m,n\ge 1$. Suppose further that $w=l+m+n$ is an {\em odd} number.
Then
\begin{equation*}
H(l,m,n;p-1)\equiv I(l,m,n)-I(n,m,l) \pmod{p}
\end{equation*}
where $I$ is defined as follows. Let $w'=w-(p-1)$ if $p<w<2p$ and $w'=w$ otherwise. Then
\begin{equation*}
I(l,m,n)=
\begin{cases}
0        &\text{ if }w\ge 2p,\text{ or if $l+m<p$ and $p<w<2p-1$},\\
1/2n  &\text{ if }w=p, 2p-1,\\
\displaystyle{ (-1)^{n+1}{w'\choose n} \frac{B_{p-w'}}{2w'}} & \text{otherwise}.
\end{cases}
\end{equation*}
In particular, if a prime $p>l+m+n$ for positive integers $l,m,n$
such that $w=l+m+n$ is odd then
\begin{equation*}
H(l,m,n;p-1)\equiv -S(l,m,n;p-1)\equiv  (-1)^{n}\left[{w \choose
l} - {w \choose n}\right] \frac{B_{p-w}}{2w}  \pmod{p}.
\end{equation*}
\end{thm}
\begin{rem} In \cite[Thm.~6.2]{H2} Hoffman gives a easier proof
than ours.
\end{rem}
\begin{proof} The last congruence follows easily from relations
\eqref{equ:rev} and \eqref{equ:partS} so we only consider the
first part. Throughout the proof all congruences are modulo $p$.
We leave the trivial cases to the interested readers and assume in
the rest of the proof that $l,m,n\ge 1$ and $w=l+m+n$ is odd.
Similar to the proof of Thm.~\ref{thm:length2} we put $L=p-1-l$,
$M=p-1-m$, and $N=p-1-n$. Then $1\le L,M,N\le p-2$. By Fermat's
Little Theorem
$$H(s_1,s_2,s_3;p-1)\equiv H(l,m,n;p-1)\equiv
\sum_{1\le k_1<k_2<k_3\le p-1} k_1^Lk_2^Mk_3^N .$$ Define the
formal power series in three variables
\begin{align*}
f(x,y,z)=&  \sum_{ r,s,t=0}^{\infty}
\left(\sum_{1\le k_1<k_2<k_3\le p-1} k_1^r  k_2^s k_3^t\right) \frac{x^r
y^sz^t}{r!s!t!}=\sum_{1\le k_1<k_2<k_3\le p-1} e^{k_1x+k_2y+k_3z} \\
=&   \frac{e^{p(x+y+z)}-e^{x+y+z}}{(e^x-1)(e^{x+y}-1)(e^{x+y+z}-1)}
-\frac{e^{x+y}(e^{pz}-e^z)}{(e^x-1)(e^{x+y}-1)(e^z-1)}\\
 &- \frac{e^x(e^{p(y+z)}-e^{y+z})}{(e^x-1)(e^y-1)(e^{y+z}-1)}
+ \frac{e^xe^y(e^{pz}-e^z)}{(e^x-1)(e^y-1)(e^z-1)} \\
=&\sum_{i,j=-1}^{\infty}  \sum_{k=0}^{\infty}
 \frac{B_{i+1}B_{j+1}(B_{k+1}(p)- B_{k+1}(1))}{(i+1)!(j+1)!(k+1)!}
  x^i(x+y)^j(x+y+z)^k\\
-&\sum_{i,j=-1}^{\infty}  \sum_{k=0}^{\infty}
 \frac{B_{i+1} B_{j+1}(1)(B_{k+1}(p)- B_{k+1}(1))}{(i+1)!(j+1)!(k+1)!}
  x^i(x+y)^jz^k \\
-&\sum_{i,j=-1}^{\infty}  \sum_{k=0}^{\infty}
 \frac{B_{i+1}(1)B_{j+1}(B_{k+1}(p)- B_{k+1}(1))}{(i+1)!(j+1)!(k+1)!}
  x^iy^j(y+z)^k\\
+&\sum_{i,j=-1}^{\infty}  \sum_{k=0}^{\infty}
 \frac{B_{i+1}(1)B_{j+1}(1)(B_{k+1}(p)- B_{k+1}(1))}{(i+1)!(j+1)!(k+1)!}
  x^iy^jz^k.
\end{align*}
We are interested in the $x^Ly^Mz^N$-term of the above sums where
$L,M,N\ge 1$. With this in mind we can safely replace all the $B_m(1)$ by
plain $B_m$. It follows that $H(l,m,n;p-1) \pmod{p}$
is the coefficient of $x^Ly^Mz^N$-term of the function $L!M!N!g(x,y,z)$ where
\begin{equation}\label{gxyz}
g(x,y,z)= \sum_{i,j=-1}^{\infty}  \sum_{k=0}^{\infty}
\frac{pB_{i+1}B_{j+1}B_k }{(i+1)!(j+1)!k!}
\Bigl( x^i(x+y)^j\big[(x+y+z)^k-  z^k\big]-x^i y^j\big[(y+z)^k -z^k\big]\Bigr).
\end{equation}
As $L+M+N=3(p-1)-w$ is odd, we know that for a term in \eqref{gxyz} to
provide a nontrivial contribution to $x^Ly^Mz^N$ it is necessary that
either $i=0$, or $j=0$, or $k=1$ because the only nonzero Bernoulli
number with odd index is $B_1=-1/2$. But $k=1$ is not an option because
all terms corresponding to $x^ay^bz$ cancel out. So we're left with only
two cases (I) $i=0$ and (J) $j=0$. We now handle them separately.

(I) When $i=0$ only the following terms really matter:
\begin{align}  \nonumber
\ & \sum_{k\ge 2} \frac{pB_1B_k}{k!} \frac{(x+y+z)^k - z^k}{x+y}  +
\sum_{j\ge 1,k\ge 2} \frac{pB_1B_{j+1}B_k}{(j+1)!k!}
  (x+y)^j\big[(x+y+z)^k - z^k \big]  \\
=&\sum_{1\le a< r\le k} \frac{pB_1B_k\cdot
     x^{a}y^{r-1-a} z^{k-r} }{r\cdot a!(r-1-a)!(k-r)!}+
\sum_{j\ge 1,k\ge 2}\sum_{c=0}^j \sum_{a+b=1}^k
\frac{pB_1B_{j+1}B_k\cdot x^{a+c}y^{j+b-c}z^{k-a-b}}
    {(j+1) c!(j-c)!a!b!(k-a-b)!} \label{1stsum}
\end{align}
where the first sum of \eqref{1stsum} comes from setting $j=-1$.

($I_0$) In the first sum of \eqref{1stsum} putting $L=a$, $M=r-1-a$,
and $N=k-r$. Then we obtain the unique term (after multiplying $L!M!N!$)
$$ \frac{pB_1B_{L+M+N+1}}{L+M+1}=\frac{pB_1B_{3(p-1)-w+1}}{2p-1-w+n}.$$
It is an easy matter to see that the denominator is a $p$-unit unless $l+m=p-1$
and the numerator is in $p\Z_p$ except when $w=2p-1$ or $w=p$.
So we find the contribution to the coefficient of $x^Ly^Mz^N$-term in $L!M!N!g(x,y,z)$
from this case:
\begin{equation*}
I_0:=\frac{pB_1B_{3(p-1)-w+1}}{2p-1-w+n}=
\begin{cases}
0  & \quad\text{ if }w\ne p, 2p-1, \\
\displaystyle{\frac{-B_1}{n}}  & \quad\text{ if }w=2p-1,\\
\displaystyle{\frac{-B_1}{p-1+n}} & \quad\text{ if }w=p,\\
\displaystyle{B_1B_{2p-1-n}} & \quad\text{ if }l+m=p-1.
\end{cases}
\end{equation*}

Now let's turn to the second sum in \eqref{1stsum}. Setting $L=a+c$, $M=j+b-c$,
and $N=k-a-b$ we see that the contribution to the coefficient of $x^Ly^Mz^N$-term is
$$I':=\sum_{r=1}^{L+M+1} \sum_{a+b=r}
\frac{pB_1 B_{L+M+1-r} B_{N+r} L!M!}{ (L+M+1-r)\cdot (L-a)!a!(M-b)!b! }.$$
Simple combinatorial argument shows that
$$\sum_{a+b=r} \frac{L!M!}{ (L-a)!a!(M-b)!b! }=
\sum_{a+b=r} {L\choose a}{M\choose b}={L+M\choose r} .$$
So we get
$$I'=\sum_{r= 1}^{L+M+1} \frac{pB_1 B_{L+M+1-r}B_{N+r} }{L+M+1-r} {L+M\choose r}.$$
If $N+L+M+1\le p-2$, i.e., $w\ge 2p$, then all terms in $I'$ are in $p\Z_p$
and we have $I'\equiv 0\pmod {p}$. Suppose $w<2p$.
It's obvious that $N+r\le L+M+N+1<3(p-1)$ and $L+M+1-r<2p-1-m-n<2p-2$.
Consequently all terms in $I'$ are in $p\Z_p$
except when ($I_1$) $L+M+1-r=p-1$, or ($I_2$) $N+r=p-1$, or ($I_3$) $N+r=2(p-1)$.
Note that $I_1$ and $I_2$ can occur at the same time if and only if $w=p$.

($I_1$) Suppose $w<2p$ and $L+M+1-r=p-1$, i.e., $r=p-l-m\ge 1$.
This occurs if and only if the followings are satisfied:
$l+m\le p-1$ and $N+r=2p-2-w\ge 1$ (so automatically $w<2p$). Then the
corresponding term in $I$ is
$$I_1 =B_1  B_{L+M+N+2-p}  {L+M\choose p-2}
=B_1  B_{2p-1-w}  {2p-2-l-m \choose p-2} .$$ However, $I_1\equiv
0$ unless $l+m=p-1$ or $w=p$. Hence
$$I_1\equiv
\begin{cases}
-B_1 B_{p-n} & \quad\text{ if $l+m=p-1$ and $n\ge 2$},\\
-B_1 B_{p-1} -B_1 & \quad\text{ if $l+m=p-1$ and $n=1$},\\
\displaystyle{\frac{B_1}{n(n-1)}} & \quad\text{ if $w=p$ and $n\ge
2$}.
\end{cases}
$$

($I_2$)  Suppose $w<2p$ and $N+r=p-1$. Because $L+M\ge r$ we see that
$w\le 2p-2$. Then the contribution is
\begin{align*}
I_2\equiv& \frac{-B_1  B_{L+M+N+2-p} }{ L+M+N+2-p}   {L+M\choose p-1-N}
\equiv \frac{-B_1B_{2p-1-w}}{2p-1-w} {2p-2-l-m\choose n} ,\\
\equiv&\begin{cases}
\text{see case ($I_1$)} & \phantom{} \quad\text{ if } w=p,\\
0 &    \quad\text{ if } l+m<p-1, \text{ or if }p<w<2p,\\
\displaystyle{\frac{B_1B_{p-n}}{p-n}}&
     \quad\text{ if $l+m=p-1$ and $n\ge 2$},\\
\displaystyle{(-1)^{n+1} \frac{B_1B_{2p-1-w}}{2p-1-w} {w+1\choose n}} &
     \quad\text{ if $l+m\ge p$ and $p<w<2p-1$, or if $w<p$}.
\end{cases}
\end{align*}

($I_3$)  Suppose $w<2p$ and $N+r=2(p-1)$. This occurs if and only if $w<p$.
Then the the corresponding term in $I$ is
\begin{equation*}
I_3\equiv\frac{-B_1  B_{L+M+N+2-p}}{L+M+N+2-p}  {L+M\choose r}
\equiv(-1)^{n}\frac{ B_1B_{p-w}}{p-w} {w\choose n-1}
\end{equation*}
by Kummer congruence.

Putting $I_0$, $I_1$, $I_2$ and $I_3$ together we get
\begin{equation}\label{equ:A}
I\equiv \begin{cases}
0   & \quad \text{ if } w\ge 2p, \\
0   & \quad \text{ if  $p<w<2p-1$ and $l+m<p$}, \\
\displaystyle{(-1)^{n+1} {w+1\choose n}} \frac{
B_{2p-1-w}}{2(w+1)} &
    \quad  \text{ if  $l+m\ge p$ and $w<2p-1$},\\
\displaystyle{1/2n} &
    \quad  \text{ if $w=p$ or $w=2p-1$},\\
\displaystyle{\frac{(-1)^{n+1}}{2w}{w\choose n} B_{p-w}} &
    \quad \text{ if } w<p.\end{cases}
\end{equation}
We notice the miracle that when $l+m=p-1$ and $n\ge 2$, all the contributions
from $I_0$, $I_1$ and $I_2$ cancel out because of Kummer congruences!

(J) When $j=0$ in \eqref{gxyz} only the following terms really matter:
\begin{align}
\ & \sum_{k\ge 2} \frac{pB_1B_k}{k!}\frac{(x+y+z)^k - (y+z)^k }{x}
+  \sum_{i\ge 1,k\ge 2} \frac{pB_1B_{i+1}B_k}{(i+1)! k!}
x^i\left((x+y+z)^k - (y+z)^k \right) \nonumber \\
=&\sum_{1\le b\le a\le k }
\frac{pB_1B_k  x^{a-1} y^b z^{k-a-b}  }{ a!b!(k-a-b)!} +
 \sum_{j\ge 1,k\ge 2} \  \sum_{a+b=1, a\ge 1}^k
\frac{pB_1B_{i+1}B_k}{(i+1)!a!b!(k-a-b)!} x^{i+a} y^b z^{k-a-b} \label{B1stsum}
\end{align}

($J_0$) In the first sum of \eqref{B1stsum} put $L=a-1$, $M=b$, and $N=k-a-b$.
After multiplying $L!M!N!$ we find the coefficient of $x^Ly^Mz^N$ is
$$J_0:= \frac{pB_1B_{3(p-1)-w+1}}{p-l}=
\begin{cases}
0  & \quad\text{ if }w\ne p, 2p-1, \\
\displaystyle{\frac{-1}{2l}}  & \quad\text{ if $w=p, 2p-1$}.
\end{cases}
$$

Let's turn to second sum in \eqref{B1stsum}. Setting $L=i+a$, $M=b$, and $N=k-a-b$
we see that the contribution from this sum is
$$J':= \sum_{a\ge 1}
\frac{pB_1 B_{L+1-a} B_{M+N+a} L!}{(L+1-a)! a! }.$$
If $L+M+N+1\le p-2$, i.e., $w\ge 2p$, then all terms are
in $p\Z_p$ and therefore $J'\equiv 0\pmod {p}$.
Suppose $w<2p$. Clearly $L+1-a\le p-2$ and
$M+N+a\le M+N+L+1<3(p-1)$.  We know that all terms in $J'$ are
in $p\Z_p$  except when ($J_1$) $M+N+a=p-1$, or ($J_2$) $M+N+a=2(p-1)$.

($J_1$) Suppose $w<2p$ and $M+N+a=p-1$. This can occur if and only if $M+N\le p-2$,
i.e., $m+n\ge p$, and $L+1-a\ge 1$ which is equivalent to $w\le 2p-2$.
If this is the case, then the contribution is
\begin{equation*}
J_1\equiv\frac{-B_1B_{2p-1-w} (p-1-l)!}{(2p-1-w)!(m+n-p+1)!}
\equiv (-1)^l\frac{B_1B_{2p-1-w}}{2p-1-w} {w+1\choose l } .
\end{equation*}

($J_2$) Suppose $w<2p$ and $M+N+a=2(p-1)$. This occurs if and only if
$L+1-a=L+M+N+3-2p=p-w\ge 1$, i.e., $w<p$. In this case
$a=m+n$ and the corresponding term in $B$ is
$$J_2 \equiv\frac{-B_1B_{p-w} (p-1-l)!}{(p-w)!(m+n)!}
  \equiv (-1)^l\frac{B_1B_{p-w}}{p-w} {w\choose l}.
$$
Combining cases ($J_0$), ($J_1$) and ($J_2$) we have
\begin{equation*}
J\equiv\begin{cases}
0   & \quad \text{ if $w\ge 2p$}, \\
0   & \quad \text{ if  $p<w<2p$ and $m+n< p$}, \\
\displaystyle{(-1)^l{w+1\choose l} } \frac{B_{2p-1-w}}{2(w+1)} &
    \quad  \text{ if $w<2p-1$ and $m+n\ge p$},\\
\displaystyle{\frac{-1}{2l}} &   \quad\text{ if $w=p, 2p-1$},\\
\displaystyle{\frac{(-1)^l}{2w} {w\choose l} B_{p-w}} &
    \quad \text{ if } w<p.
\end{cases}
\end{equation*}

The theorem now follows from an easy simplification process.
\end{proof}

In the case when $w=r+s+t$ is even the expression is much more
complicated. For future reference we provide the following
computation when the prime $p>w+1$. The method is similar to that
used by Hoffman to compute the length 2 case in \cite{H2}. Recall
that by Bernoulli polynomials we have
\begin{equation}\label{equ:sump}
\sum_{i=1}^n i^d=\sum_{a=0}^d {d+1\choose a}\frac{B_a}{d+1}
n^{d+1-a}.
\end{equation}
So modulo $p$ we have by Fermat's Little Theorem
\begin{align*}
    H(r,s,t;p-1)\equiv&\sum_{i=1}^{p-1}\frac{1}{i^t}\sum_{j=1}^{i-1}
    j^{p-1-s} \sum_{k=1}^{j-1} k^{p-1-r}\\
    \equiv&\sum_{i=1}^{p-1}\frac{1}{i^t}\sum_{j=1}^{i-1}
     \sum_{a=0}^{p-1-r}  {p-r\choose a}\frac{B_a }{p-r} j^{\gk(a)+p-r-s-a}\\
    \equiv&\sum_{i=1}^{p-1}\frac{1}{i^t}
     \sum_{a=0}^{p-1-r}  {p-r\choose a} \frac{B_a}{p-r}\sum_{j=1}^{i-1} j^{\gk(a)+p-r-s-a}\\
    \equiv&  \sum_{i=1}^{p-1}\sum_{a=0}^{p-1-r} {p-r\choose a} \frac{B_a}{p-r}\cdot\\
    \ &  \sum_{b=0}^{\gk(a)+p+1-r-s-a}  {\gk(a)+p+1-r-s-a\choose
    b}\frac{B_b\cdot
      i^{\gk(a)+p+1-w-a-b}}{\gk(a)+p+1-r-s-a},
\end{align*}
where $\gk(a)=0$ if $a\le p-r-s$ and $\gk(a)=p-1$ if $a>p-r-s$.
Now take the sum of powers of $i$ first and observe that
$\sum_{i=1}^{p-1} i^l\equiv 0 \pmod{p}$ unless $l\equiv
0\pmod{p-1}$. Note that it is impossible to have
$\gk(a)+p+1-w-a-b\equiv 0 \pmod{p-1}$ if $a>p+1-w$, unless
$a>p-r-s$. Hence we get:
\begin{align}\notag
H(r,s,t;p-1)& \equiv -\sum_{a=0}^{p+1-w}  {p-r\choose
    a}\frac{B_a}{p-r}  {p+1-r-s-a\choose t}
    \frac{B_{p+1-w-a}}{p+1-r-s-a}\\
    &\quad -\sum_{a=\max\{p+1-r-s,p+2-w\}}^{p-1-r}  {p-r\choose
    a}\frac{B_a}{p-r}  {2p-r-s-a\choose t} \frac{B_{2p-w-a}}{2p-r-s-a}\notag\\
& \aligned \equiv &-\sum_{a=0}^{p+1-w}  (-1)^{a+t}  {r+a \choose
    a}\frac{B_a}{r+a}  {w+a-1\choose t}
    \frac{B_{p+1-w-a}}{w+a-1}\\
   &-\sum_{a=\max\{p+1-r-s,p+2-w\}}^{p-1-r} (-1)^{a+t}  {r+a \choose
    a}\frac{B_a}{r+a}  {w+a\choose t} \frac{B_{2p-w-a}}{w+a}.\endaligned
 \label{equ:Hrstgeneral}
\end{align}
When $w$ is odd we recover Thm.~\ref{thm:oddlength3} when $p>w+1$
by noticing that there are only two nontrivial terms in
\eqref{equ:Hrstgeneral} corresponding to $a=1$ and $a=p-w$.

\subsection{Some remarkable cases of length 3 with even weight}\label{sec:d3wt4}
Applying \eqref{equ:s1s2} to $(s_1, s_2)=(1,3)$ together with
shuffle product $H(1,1;p-1)\cdot H(2;p-1)$ we get
\begin{equation}\label{equ:cyclic}
H(2,1,1;p-1)+H(1,2,1;p-1)+H(1,1,2;p-1) \equiv \frac{4p}{5}B_{p-5}
\pmod{p^2}.
\end{equation}
By Lemma~\ref{lem:revorder} we can even see that
 \begin{equation}\label{equ:112=211}
H(2,1,1;p-1)\equiv H(1,1,2;p-1) \pmod{p}.
\end{equation}However, is it true that
in fact all of these sums are congruent to 0 mod $p$? Now from
\eqref{equ:Hrstgeneral} we can compute easily that modulo $p$
\begin{align}
 H(1,2,1;p-1)\equiv&\sum_{a=0}^{p-3} B_aB_{p-3-a}, \label{equ:H1211}\\
H(1,1,2;p-1)\equiv&-\sum_{a=0}^{p-3} \frac{2+a}{2} B_aB_{p-3-a}, \label{equ:H1212}\\
H(2,1,1;p-1)\equiv&\sum_{a=0}^{p-3} \frac{1+a}{2}
 B_aB_{p-3-a}. \label{equ:H1213}
\end{align}
Set
$$A:=\sum_{a=0}^{p-3} B_aB_{p-3-a},\quad B :=\sum_{a=0}^{p-3} a
 B_aB_{p-3-a}.$$
Then from \eqref{equ:cyclic}, \eqref{equ:H1211} to
\eqref{equ:H1213}
\begin{equation*}
 A/2 \equiv 0\pmod{p}.
\end{equation*}
Moreover, from \eqref{equ:112=211}, \eqref{equ:H1212} and
\eqref{equ:H1213}
$$-(A+B/2) \equiv A/2+B/2\pmod{p}.$$
Consequently we have
\begin{cor} \label{equ:A=B=0}
For every prime $p\ge 7$ we get
$$\sum_{a=0}^{p-3} B_aB_{p-3-a}\equiv \sum_{a=0}^{p-3} a
 B_aB_{p-3-a}\equiv 0 \pmod{p}.$$
Therefore we have
\begin{equation}\label{equ:112211=0}
H(1,2,1;p-1)\equiv H(1,1,2;p-1)\equiv H(2,1,1;p-1)\equiv
0\pmod{p}.
\end{equation}
\end{cor}

Using Maple we further find the following very stimulating example
$$H(1,2,1;36)=\frac{2234416196881673576349577192603}
{1151149136943530805554073600000}\equiv 0\pmod{37^2}.$$
\begin{prop}
For all prime $p\ge 7$ we have
\begin{align}
H(1,2,1;p-1)\equiv &-\frac{9}{10}pB_{p-5}\pmod{p^2},\label{cong:H121}\\
H(2,1,1;p-1)\equiv &\phantom{-} \frac{3}{5}pB_{p-5} \pmod{p^2}, \label{cong:H211}\\
H(1,1,2;p-1)\equiv &\phantom{-} \frac{11}{10}pB_{p-5} \pmod{p^2}.
\label{cong:H112}
\end{align}
And
\begin{align}
S(1,2,1;p-1)\equiv &-\frac{9}{10}pB_{p-5}\pmod{p^2},\label{cong:S121}\\
S(2,1,1;p-1)\equiv &\phantom{-} \frac{11}{5}pB_{p-5} \pmod{p^2}, \label{cong:S211}\\
S(1,1,2;p-1)\equiv &\phantom{-} \frac{3}{10}pB_{p-5} \pmod{p^2}.
\label{cong:S112}
\end{align}
\end{prop}
\begin{proof} Omitting $H(\cdots; p-1)$ we let $A=H(1,2,1),
B=H(2,1,1)$ and $C=H(1,2,2)$. Equation \eqref{equ:cyclic} says
that
\begin{equation} \label{equ:cyclic1}
A+B+C \equiv \frac{4}{5}pB_{p-5} \pmod{p^2}.
\end{equation}
Now by the shuffle relations
\begin{equation}\label{equ:H1H21}
H(1)H(2,1)=A+2B+H(3,1)+H(2,2),\\
\end{equation}
 From Thm.~\ref{thm:nonhomol=2} we get
\begin{equation}\label{equ:A+2B}
 A+2B \equiv   \frac{3}{10}p B_{p-5} \pmod{p^2}\\
\end{equation}
Hence it suffice to show \eqref{cong:H121}. From \eqref{equ:SH} we
see that
\begin{alignat}{2}
S(1,2,1)=&\,H(1,2,1)+H(3,1)+H(1,3)+H(4) \notag\\
& \hskip1cm =H(1,2,1)+H(1)\cdot H(3)\equiv H(1,2,1) & \pmod{p^2}\ \,  \label{S121sim}\\
 S(2,2)=&\, H(2,2)+H(4)=-H(2,2)+H(2)^2\equiv -H(2,2) & \pmod{p^2}.  \label{S22sim}
\end{alignat}
Now because $(1,2,1)^\ast=(2,2)$ by Thm.~\ref{cong:modp2} we have
\begin{equation}\label{cong:22S121}
-S(2,2)\equiv
S(1,2,1)+p\Big(H(1,2,1,1)+H(3,1,1)+H(1,3,1)+H(4,1)\Big)
\pmod{p^2}.
\end{equation}
By Cor.~\ref{cor:palin} $H(1,3,1)\equiv 0 \pmod{p}.$ So using
expressions \eqref{S121sim} and \eqref{S22sim} we can simplify the
preceding congruence to
\begin{equation}\label{cong:121}
H(2,2)\equiv H(1,2,1)+p\Big(H(1,2,1,1)+H(3,1,1)+H(4,1)\Big)
\pmod{p^2}.
\end{equation}
Now that $(1,2,1,1)^\ast=(2,3)$ we find from congruence
\eqref{equ:conv} that modulo $p$
\begin{multline*}
 -S(2,3)\equiv
S(1,2,1,1)\equiv H(1,2,1,1)+H(3,1,1)+H(1,3,1)\\
 +H(1,2,2)+H(4,1)+H(3,2)+H(1,4)+H(5) \pmod{p}.
 \end{multline*}
Note that
$$H(1,3,1)\equiv 0,H(4,1)+H(1,4)=H(1)\cdot H(4)- H(5) \equiv 0
\pmod{p}.$$
 So \eqref{equ:SH} implies that
 $$
 -H(2,3)\equiv  -S(2,3)\equiv
 H(1,2,1,1)+H(3,1,1)
 +H(1,2,2)+H(3,2)  \pmod{p}.$$
Namely
$$\aligned
H(1,2,1,1) \equiv& -H(2,3)-H(3,2)-H(3,1,1)
 -H(1,2,2) & \pmod{p}\  \\
 \equiv &\, H(5)-  H(3,1,1) -H(1,2,2) &   \pmod{p} \  \\
 \equiv & -H(3,1,1) -H(1,2,2)  &\pmod{p}.
\endaligned$$
Plugging this into \eqref{cong:121} we see that
$$H(2,2)\equiv
 H(1,2,1)+p\Big(H(4,1)-H(1,2,2) \Big) \pmod{p^2}.
$$
Using Thm.~\ref{thm:zc} and Thm.~\ref{thm:nonhomol=2} we can now
compute easily that
$$H(2,2)\equiv -\frac25 p B_{p-5} \pmod{p^2},\quad
 H(4,1)\equiv -B_{p-5}, \quad H(1,2,2)\equiv
 -\frac32 B_{p-5} \pmod{p}.$$
These lead to congruence \eqref{cong:H121}. Then congruence
\eqref{cong:S121} follows from \eqref{S121sim}.

We now can solve \eqref{equ:cyclic1} and \eqref{equ:A+2B} to get
congruences \eqref{cong:H211} and \eqref{cong:H112}. Finally,
\eqref{equ:partS} yields
\begin{align*}
    S(2,1,1)= &\,H(1,1,2)-H(1)H(1,2)-H(1,1) H(2)+H(1)^2H(2)\equiv \frac{11}{10}pB_{p-5} \pmod{p^2},\\
    S(1,1,2)= &\,H(2,1,1)-H(2)H(1,1)-H(2,1)H(1)+H(2)H(1)^2 \equiv  \frac{3}{5}pB_{p-5}\pmod{p^2}.
\end{align*}
We have completed the proof of the proposition.
\end{proof}

By going through the proof of Thm.~\ref{thm:oddlength3} or the
proof leading to \eqref{equ:Hrstgeneral} we can obtain the
following result. We leave the details of the proof of it to the
interested readers.

\begin{prop}\label{equ:435}
For all prime $p\ge 17$ we have
\begin{equation*}
H(4,3,5;p-1)\equiv  H(5,3,4;p-1)\equiv 0\pmod{p}.
\end{equation*}
\end{prop}

In \cite{pmod}, by using the shuffle relations and Hoffman's
convolution we will study the mod $p$ structure of the multiple
harmonic sums for lower weights. In particular, we will prove
Prop.~\ref{equ:435} and congruences like
 $$S(2,3,2,3,2;p-1)\equiv S(2,3,3,2,2;p-1)\equiv  0\pmod{p}.$$

\begin{prob} \label{prob:H435}
Numerical evidence shows that if $(r,s,t)\ne (5,3,4), (4,3,5)$ and
$r+s+t\ge 6$ is even then the density of primes $p$ such that
$p\nmid H(r,s,t;p-1)$ among all primes is always 1; however, there
is always $p$ such that $p|H(r,s,t;p-1)$. Can  one generalize the
formula in Thm.~\ref{thm:oddlength3} to prove this?
\end{prob}

\subsection{Some congruences modulo prime squares}\label{sec:congBer}
We know from Thm.~\ref{thm:oddlength3} that $H(r,s,r;p-1)\equiv 0
\pmod{p}$ if $s$ is an odd number and $p>2r+s$. Modulo $p^2$ we
have (using substitution of indices $k\to p-k$)
\begin{equation*}
   H(r,s,r;p-1)\equiv
   -H(r,s,r;p-1)-p\big[2rH(r+1,s,r;p-1)+sH(r,s+1,r;p-1)\big].
\end{equation*}
Hence
\begin{equation}\label{equ:modpsq}
   H(r,s,r;p-1)\equiv
   -p\big[ rH(r+1,s,r;p-1)+\frac{s}{2}H(r,s+1,r;p-1)\big]\pmod{p^2}.
\end{equation}
\begin{prop} For all prime $p>5$ we have
\begin{align}\label{equ:H131}
   H(1,3,1;p-1)\equiv &\, S(1,3,1;p-1)\equiv 0 \quad \pmod{p^2}\\
\label{equ:H212}
   H(2,1,2;p-1)\equiv &\, S(2,1,2;p-1)\equiv -\frac{1}{3}B_{p-3}^2 p \pmod{p^2}.
\end{align}
For all prime $p>7$ set $b_8(p) =(5 S(6,1,1;p-1) +
B_{p-5}B_{p-3})/2$ then we have
\begin{align}\label{equ:H151}
   H(1,5,1;p-1)\equiv  &\, S(1,5,1;p-1)\equiv b_8(p)p
   \pmod{p^2}\\
   H(2,3,2;p-1)\equiv  &\,S(2,3,2;p-1)\equiv
   \, 4b_8(p)p \pmod{p^2},\label{equ:H232}\\
   H(3,1,3;p-1)\equiv  &\, S(3,1,3;p-1)\equiv b_8(p) p
   \pmod{p^2}. \label{equ:H313}
\end{align}
\end{prop}
\begin{proof} Notice that
$$S(r,s,r;p-1)\equiv H(r,s,r;p-1)+H(r+s)H(r) \equiv
H(r,s,r;p-1)\pmod{p}.$$
 From \eqref{equ:modpsq} we have
$$\aligned
   2H(1,3,1;p-1)\equiv  &\,-p\big [2H(2,3,1;p-1)+3H(1,4,1;p-1)\big] \pmod{p^2},\\
   2H(2,1,2;p-1)\equiv  &\,-p\big [4H(3,1,2;p-1)+H(2,2,2;p-1)\big] \pmod{p^2},\\
   2H(1,5,1;p-1)\equiv  &\,-p\big [2H(2,5,1;p-1)+5H(1,6,1;p-1)\big] \pmod{p^2},\\
  2H(2,3,2;p-1)\equiv   &\,-p\big [4H(3,3,2;p-1)+3H(2,4,2;p-1)\big] \pmod{p^2},\\
   2H(3,1,3;p-1)\equiv  &\,-p\big [6H(4,1,3;p-1)+H(3,2,3;p-1)\big] \pmod{p^2}.
\endaligned$$
By \cite[Thm.~7.2]{H2} we know that
\begin{align*}
S(2,3,1;p-1)\equiv &\,3S(4,1,1;p-1) &\,\pmod{p} ,\\
S(1,4,1;p-1)\equiv &\,-2S(4,1,1;p-1) &\,\pmod{p},\\
S(3,1,2;p-1)\equiv &\,-S(4,1,1;p-1)\equiv
\frac{1}{6}B_{p-3}^2&\,\pmod{p}
\end{align*}
which yield \eqref{equ:H131} and \eqref{equ:H212}. Similarly,
 \eqref{equ:H151} follows from \cite[Thm.~7.4]{H2} because
 \begin{align*}
2S(2,5,1;p-1)\equiv &\,-2S(5,2,1;p-1)-2S(5,1,2;p-1) &\, \pmod{p},\\
 \equiv &\, 5 S(6,1,1;p-1)-B_{p-5}B_{p-3} &\, \pmod{p},\\
S(1,6,1;p-1)\equiv &\,-2 S(6,1,1;p-1) &\, \pmod{p} .
 \end{align*}
Also from \cite[Thm.~7.4]{H2}
 \begin{align*}
S(3,3,2;p-1)\equiv &\, 4b_8(p) &\, \pmod{p},\\
S(2,4,2;p-1)\equiv &\,-2 S(4,2,2;p-1) \equiv -2S(3,3,2;p-1) &\, \pmod{p}, \\
S(3,2,3;p-1)\equiv &\, -2S(3,3,2;p-1) &\, \pmod{p},\\
S(4,1,3;p-1)\equiv &\,  b_8(p) &\,\pmod{p} .
 \end{align*}
These lead to the last two congruence of the proposition
immmediately.
\end{proof}

By similar argument we can compute $H(r,s,r) \pmod{p^2}$ for odd
$s$ if we know the values $S(\ors;p-1) \pmod{p}$ for all $\ors$ of
length 3 and weight $2r+s+1$. When we apply this argument to
$\ors=(1,3,1,3)$ we find
\begin{prop} The following two congruences are equivalent:

(i) For all primes $p>8$
$$H(1,3,1,3; p-1) \equiv -\frac{31}{72}
pB_{p-9} \pmod{p^2},$$

(ii)  For all primes $p>8$
$$S(6,1,1,1;p-1)\equiv -\frac{1}{54}B_{p-3}^3 -\frac{1889}{648} B_{p-9} \pmod{p}.$$
\end{prop}
We have verified the congruences in the proposition for all primes
$p$ such that $10<p<2000$.

\subsection{Some congruences of Bernoulli numbers}\label{sec:congBer2}
 From Cor.~\ref{equ:A=B=0} we see that for every prime $p\ge 7$ we
 have
$$\sum_{a=0}^{p-3} B_aB_{p-3-a}\equiv \sum_{a=0}^{p-3} a
 B_aB_{p-3-a}\equiv 0 \pmod{p}.$$
Can we generalize this? The answer turns out to be affirmative.
\begin{prop} \label{prop:berp-3}
For every prime $p\ge 9$ we have
$$\sum_{a=0}^{p-5} B_aB_{p-5-a}\equiv  -\frac23 B_{p-3}^2 \pmod{p}.$$
\end{prop}
\begin{proof}
By \eqref{equ:Hrstgeneral} we have for any even number $n$
\begin{equation}\label{equ:h141}
    H(1,n,1;p-1)\equiv \sum_{a=0}^{p-n-1}(-1)^a
    B_aB_{p-n-1-a}+\sum_{a=p-n}^{p-2} (-1)^a
    B_aB_{2p-n-2-a}
\end{equation}
Taking $n=4$ and comparing with \cite[Thm.~7.2]{H2} we get
$$\sum_{a=0}^{p-5} B_aB_{p-5-a}+ B_{p-3}^2\equiv
 H(1,4,1)\equiv S(1,4,1)\equiv \frac13 B_{p-3}^2 \pmod{p}.$$
This proves the proposition.
\end{proof}

The following result is straight-forward.
\begin{prop} For all positive number $n$ and prime $p>n+3$ we have
\begin{multline*}
 H(1,1,1,n;p-1) \equiv -(-1)^n H(n,1,1,1;p-1) \equiv (-1)^n S(1,1,1,n;p-1)\equiv -S(n,1,1,1;p-1) \\
 \equiv \sum_{a=0}^{p-2}\sum_{b=0}^{p-n-1} (-1)^{b+n} {
 a+b\choose b} {a+b+n\choose n} \frac{B_aB_bB_{ p-n-1-a-
b}}{(a+1)(a+b+1)} \pmod{p}.
\end{multline*}
\end{prop}

\subsection{Multiple harmonic sums of arbitrary length}\label{sec:conv}
To prove the main result in this section let us recall the
Bernoulli polynomial $B_m(x)$ which is defined by the following
generating function
$$\frac{te^{xt}}{e^t-1}=\sum_{m=0}^\infty B_m(x)\frac{t^m}{m!}.$$
These polynomials satisfy (see \cite[p.248]{IR}):
\begin{align}\label{equ:Berpoly1}
B_m(x)=&\sum_{k=0}^{m}  {m\choose k} B_k x^{m-k}\\
B'_m(x)=& m B_{m-1}(x) \quad \forall m\ge 1. \label{equ:Berpoly2}
\end{align}
\begin{lem} For any even integer $n\ge 0$ and prime $p\ge 3n+7$ we
have
\begin{align} \label{equ:B2=0}
\sum_{a,b\ge 0,\, a+b=p-3n-3} \frac{B_b}{p-b}{p-b\choose n+1}
 \frac{B_a}{p-a} {p-a\choose n+2}\equiv 0 & \pmod {p},\\
 \sum_{a,b\ge 0,\, a+b=p-3n-3} \frac{B_b}{p-b}{p-b\choose n+1}
 \frac{B_a}{p-a} {p-a\choose n+1} \equiv  0 &\pmod {p}.
 \label{equ:DerB2=0}
\end{align}
\end{lem}
\begin{proof} Throughout the proof all congruences are modulo $p$. First we set
\begin{align*}
A=&\sum_{a,b\ge 0,\, a+b=p-3n-3} \frac{B_b}{p-b}{p-b\choose n+1}
 \frac{B_a}{p-a} {p-a\choose n+2},\\
B=& \sum_{a,b\ge 0,\, a+b=p-3n-3} \frac{B_b}{p-b}{p-b\choose n+1}
 \frac{B_a}{p-a} {p-a\choose n+1}.
 \label{equ:DerB2=0}
\end{align*}
We have
\begin{equation}\label{equ:lhsnn+1}
A\equiv  - \sum_{a,b\ge 0,\, a+b=p-3n-3} \frac{n+1+a}{n+2} \cdot
\frac{B_b}{p-b}{p-b\choose n+1} \frac{B_a}{p-a} {p-a\choose n+1}.
\end{equation}
Now exchange the index $a$ and $b$ in $A$ we get
\begin{align}\notag
A=& \sum_{a,b\ge 0,\, a+b=p-3n-3} \frac{B_b}{p-b}{p-b\choose n+2}
 \frac{B_a}{p-a} {p-a\choose n+1} \\
\equiv& \sum_{a,b\ge 0,\, a+b=p-3n-3} \frac{2n+2+a}{n+2}\cdot
\frac{B_b}{p-b}{p-b\choose n+1}
 \frac{B_a}{p-a} {p-a\choose n+1}. \label{equ:lhsnn+1a}
\end{align}
Adding \eqref{equ:lhsnn+1} to \eqref{equ:lhsnn+1a} we find
 \begin{equation}\label{2a=b}
2A\equiv \sum_{a,b\ge 0,\, a+b=p-3n-3} \frac{n+1}{n+2}\cdot
\frac{B_b}{p-b}{p-b\choose n+1}
 \frac{B_a}{p-a} {p-a\choose n+1} \equiv   \frac{n+1}{n+2}B.
\end{equation}
Now for any nonnegative integer $m$ and $b$ such that $b+m+1<p$ we
have
\begin{equation}\label{equ:binoequiv}
 \frac{(-1)^m}{p-b}{p-b\choose
m+1} \equiv (-1)^{b}m{p-m-1\choose b}.
\end{equation}
We get
\begin{equation}\label{equ:DerB2=0a}
B\equiv n^2 \sum_{a,b\ge 0,\, a+b=p-3n-3} B_b {p-n-1\choose b}
 B_a  {p-n-1\choose a}
\end{equation}
because all the terms with odd index $a$ are zero. Let us denote
by Coeff$[i, f(x)]$ the coefficient of $x^i$ in the polynomial
$f(x)$. Then \eqref{equ:DerB2=0a} means that
 \begin{align*}
 B\equiv &\, n^2 \text{Coeff}[ p+n+1, B^2_{p-n-1}(x) ]\\
  \equiv &\frac{ n^2}{p+n+1}\text{Coeff}[ p+n, \big(B^2_{p-n-1}(x)
  \big)']\\
 \equiv &\frac{n^2}{n+1}\text{Coeff}[ p+n,
 2(p-n-1)B_{p-n-1}(x)B_{p-n-2}(x)]
\end{align*}
by \eqref{equ:Berpoly2}. From \eqref{equ:Berpoly1} this yields
\begin{align*}
B\equiv& -2n^2 \sum_{a,b\ge 0,\, a+b=p-3n-3} B_b {p-n-1\choose b}
 B_a  {p-n-2\choose a} \\
 \equiv & \frac{2n}{n+1} \sum_{a,b\ge 0,\, a+b=p-3n-3} \frac{B_b}{p-b}{p-b\choose n+1}
 \frac{B_a}{p-a} {p-a\choose n+2}
\end{align*}
by \eqref{equ:binoequiv}. Therefore
\begin{equation}\label{2b=a}
    B\equiv \frac{2n}{n+1} A.
\end{equation}
Now the lemma follows easily by \eqref{2a=b} and \eqref{2b=a}.
\end{proof}

\begin{thm} \label{thm:conv}
Let $p$ be a prime and $\ors\in \N^l$. Assume $p>|\ors|+2$. Then
$$H(\ors;p-1)\equiv S(\ors;p-1)\equiv 0\pmod{p}$$
if $\ors$ has the following forms:
\begin{enumerate}

\item \label{m2n} $\ors=\big (1^m,2,1^n\big )$ for $m,n\ge 0$ and
$m+n$ is even.

\item \label{nn-1} $\ors= \big(1^{n},2,1^{n-1},2,1^{n+1}\big )$
where $n\ge 2$ is even.

\item \label{n-1n} $\ors= \big(1^{n+1},2,1^{n-1},2,1^{n}\big)$
where $n\ge 0$.

\item \label{nnn} $\ors= \big(1^{n},2,1^n,2,1^n\big )$ where $n\ge
2$ is even.

Furthermore, in the first and last cases we also have
$S(\ors;p-1)\equiv 0\pmod{p}.$
\end{enumerate}
\end{thm}
\begin{proof} We omit $p-1$ in $S(-;p-1)$ and $H(-;p-1)$ and
assume all congruences are modulo $p$ throughout this proof.

\eqref{m2n}. Let $\ors=\big (1^m,2,1^n\big )$. For future
reference we first allow $m$ and $n$ to have different parity.
Clearly we have $\ors^\ast=(m+1,n+1)$. So
$|\ors|=|\ors^\ast|=m+n+2$, $l(\ors)=m+n+1$, and $l(\ors^\ast)=2$.
By \cite[Thm.~6.7]{H1} we have
\begin{equation}\label{*=bar}
    S(\ors^\ast)\equiv -S(\ors).
\end{equation}
Since $|\ors^\ast|-l(\ors^\ast)=m+n$ we have by \eqref{equ:partS}
\begin{equation}\label{*=0}
  (-1)^{m+n} S(\ors^\ast)\equiv (-1)^{|\ors^\ast|-l(\ors^\ast)} S(\ors^\ast)
    \equiv -H(\ors^\ast)+H(m+1)\cdot H(n+1)\equiv -H(\ors^\ast)
\end{equation}
by Thm.~\ref{thm:Bayat}. We now apply \eqref{equ:partS} to $\ors$
and get
\begin{equation}\label{equ:plain}
    -S(\ors)\equiv \sum_{\bigsqcup_{j=1}^l \ors_j=\ors}
(-1)^l\prod_{j=1}^l H(\ors_j).
\end{equation}
Now if $l\ge 2$ then one of $\ors_j=1^d$ for some positive $d$ so
that $H(\ors_j)\equiv 0.$ Thus
$$S(\ors)\equiv H(\ors).$$
Combining this with \eqref{*=bar} and \eqref{*=0} we get
\begin{equation}\label{equ:s=sast}
      H(\ors)\equiv (-1)^{m+n}  H(\ors^\ast).
\end{equation}
When $m+n$ is even we get
$$S(\ors)\equiv H(\ors) \equiv  H(\ors^\ast)  \equiv 0 .$$
by Thm.~\ref{thm:nonhomol=2} which proves case \eqref{m2n}.

\eqref{nn-1} Let $\ors= \big(1^{n},2,1^{n-1},2,1^{n+1}\big )$
where $n\ge 2$ is even. Then $\ors^\ast=(n+1,n+1,n+2).$ So
$|\ors|=|\ors^\ast|=3n+4$, $l(\ors)=3n+2$, and $l(\ors^\ast)=3$.
Since $n$ is even we have by applying \eqref{equ:partS} to
$\ors^\ast$
\begin{equation}\label{case2*=0}
    S(\ors^\ast)\equiv
    H(\ors^\ast)-H(2n+2,n+2)-H(n+1,2n+3)+H(3n+4)\equiv H(\ors^\ast)
\end{equation}
by Thm.~\ref{thm:nonhomol=2} and Thm.~\ref{thm:Bayat}. Applying
\eqref{equ:partS} to $\ors$ and using the fact that $H(1^d)\equiv
0$ for any $d$ we have
\begin{align*}
    S(\ors)\equiv (-1)^{|\ors|-l(\ors)} S(\ors)\equiv &
    -H(\ors)+\sum_{a=0}^{n-1} H(1^n,2,1^a)\cdot
    H(1^{n-1-a},2,1^{n+1})\\
    \equiv & -H(\ors)+\sum_{a=0}^{n-1} H(n+1,a+1)\cdot H(n-a,n+2)
\end{align*}
by \eqref{equ:s=sast}. Hence by \eqref{*=bar} and \eqref{case2*=0}
\begin{equation}\label{equ:Hors=nnn}
H(\ors)\equiv H(n+1,n+1,n+2)+ \sum_{a=0}^{n-1}  H(n+1,a+1)\cdot
H(n-a,n+2).
\end{equation}
We know that for all $j,k<p$ we have
$$\sum_{a=0}^{p-2} (k/j)^a \equiv
\begin{cases} 0 \pmod{p} &\text{ if }j\ne k,\\
 -1 \pmod{p} &\text{ if }j= k.
 \end{cases} $$
It follows that
\begin{align*}
 \sum_{a=0}^{p-2} H(n+1,a+1)\cdot
 H(n-a,n+2)=& \sum_{a=0}^{p-2} \sum_{1\le i<j<p} \, \sum_{1\le k<l<p}
 \frac{1}{i^{n+1}j^{a+1}} \frac{1}{k^{n-a}l^{n+2}}\\
    \equiv &-\sum_{1\le i<j=k<l<p}
 \frac{1}{i^{n+1}j^{n+1}l^{n+2}} \\
 =& -H(n+1,n+1,n+2).
\end{align*}
Together with \eqref{equ:Hors=nnn} we see that
 \begin{align*}
H(\ors)\equiv &-\sum_{a=n}^{p-2}
 H(n+1,a+1)\cdot H(n-a,n+2) \\
 \equiv &-\sum_{a=2n+1}^{p-n-2} \frac{(-1)^a B_{p-n-2-a} }{n+2+a}
 {n+2+a\choose n+1} \frac{B_{a-2n-1}}{p+2n+1-a} {p+2n+1-a\choose
 n+2}
\end{align*}
by Thm.~\ref{thm:length2}. Under substitution $a\rightarrow
2n+1+a$ we get:
\begin{equation}\label{equ:Hors=w}
H(\ors)\equiv \sum_{a=0}^{p+1-w} (-1)^a {w-1+a \choose
    n+1}\frac{B_{p+1-w-a}}{w-1+a}{p-a\choose n+2}
    \frac{B_a}{p-a}\equiv 0\pmod{p}
\end{equation}
by \eqref{equ:B2=0}. This combined with \eqref{case2*=0} completes
the proof of case \eqref{nn-1}.

\eqref{n-1n} It follows from \eqref{nn-1} by taking $\bar\ors$.

\eqref{nnn} When $n=0$ or $n$ is odd this follows from
Thm.~\ref{thm:p-1} and Lemma~\ref{lem:revorder}, respectively.
When $n\ge 2$ is even the proof is almost the same as that of case
\eqref{nn-1} except at the end one need resort to
\eqref{equ:DerB2=0}. The congruence for $S$ follows from the fact
that $\big(1^{n},2,1^n,2,1^n\big )^\ast=(n+1,n+2,n+1)$ and
therefore
$$S(\ors^\ast)\equiv H(\ors^\ast)\equiv 0$$
by Cor.~\ref{cor:palin}.
\end{proof}

\begin{rem} Note that in cases \eqref{nn-1} and \eqref{n-1n} we
usually have $S(\ors)\not\equiv 0\pmod{p}$. For example, in case
\eqref{nn-1} we have  $S(\ors) \equiv -S(\ors^\ast)\equiv
-H(\ors^\ast) \equiv H(n+1,n+1,n+2) \pmod{p}$ by \eqref{case2*=0}.
We know that $H(3,3,4;12)\equiv 8\pmod{13}, H(5,5,6;18)\equiv
15\pmod{19}$ and $H(7,7,8;28)\equiv 26 \pmod{29}$.
\end{rem}

\begin{thm}\label{thm:r+seven}
Let $\ors=\{r,s\}^n$ for some $r,s\ge 1$ and $p\ge |\ors|$ be a
prime. Then
$$H(\ors;p-1)\equiv S(\ors;p-1)\equiv 0\pmod{p}$$
if either (i) $n=1,2$, both $r$ and $s$ are even, or (ii) $n$ is
any positive integer, both $r$ and $s$ are odd.
\end{thm}
\begin{rem}  When $n=1$ this is Thm.~\ref{thm:length2}. When $n=2$,
this has been confirmed by Hoffman (see the remarks after
\cite[Thm.~6.3]{H2}).
\end{rem}
\begin{proof} By the above remark we may assume $r$ and $s$ are odd and
proceed by induction on $n$. In the following we will drop $p-1$
again. By the shuffle relations and equation \eqref{equ:SH} it is
straightforward to verify that
\begin{align}
S(\{r,s\}^n)= &  \sum_{\ort \preceq \{r,s\}^n} H(\ort) \notag\\
 =&\,H(\{r,s\}^n)+ H(r+s)\cdot \sum_{\ort \preceq \{r,s\}^{n-1}}
 H(\ort) \notag \\
 =& \,H(\{r,s\}^n)+ H(r+s)\cdot S(\{r,s\}^{n-1})\equiv
 H(\{r,s\}^n) \pmod{p}\label{equ:Srs}
\end{align}
by induction assumption. On the other hand, from \eqref{equ:partS}
we find that
$$ S(\{r,s\}^n)\equiv \sum_{
 \bigsqcup_{i=1}^l   \ors_i=\{r,s\}^n} (-1)^l H(\ors_1)\cdots
 H(\ors_l) \pmod{p}.$$
Now if $l>1$ then $\ors_1=\{r,s\}^d$ for some $d<n$ in which case
$ H(\ors_1)\equiv 0 \pmod{p}$ by induction assumption, or else,
$\ors_1=\{r,s\}^d\sqcup \{r\}$, in which case $\ors_1$ is a
palindrome of odd weight and hence $ H(\ors_1)\equiv 0 \pmod{p}$
by Cor.~\ref{cor:palin}. Consequently
$$ S(\{r,s\}^n)\equiv -H(\{r,s\}^n) \pmod{p}.$$
Together with \eqref{equ:Srs} it shows that $ H(\{r,s\}^n)\equiv
S(\ors;p-1)\equiv 0 \pmod{p}$ and the theorem is proved.
\end{proof}

When  both $r$ and $s$ are even but $n>2$ the theorem does not
hold in general. For example,
$$H(\{2,4\}^3;22)\equiv 21 \pmod{23},\quad H(\{2,4\}^4;28)\equiv 20
\pmod{29}.$$

\subsection{Some conjectures in general cases}
When the length $l\ge 4$ numerical evidence up to length 40 shows
the following conjecture is true.
\begin{conj}\label{conj:cases}
Let $\ors\in \N^l$ and $p\ge |\ors|$ be a prime. Then
$$H(\ors;p-1)\equiv 0\pmod{p}$$
if $\ors$ has one of the following forms:
\begin{enumerate}
\item $\ors=\big (\big \{1^m,2,1^n,2 \big \}^q,1^m,2,1^n \big )$
for $q,m, n\ge 0$, where either (i) $q$ is odd, or (ii) $q$ is
even and $m+n$ is even.

\item $\ors=\big (\{2\}^m,\big\{3,\{2\}^m\big \}^n\big )$ for
$m,n\ge 0$.

\item $\ors=\big (1,\{2\}^m,\big\{1,\{2\}^{m+1}\big
\}^n,1,\{2\}^m,1\big )$ for $m,n\ge 0$ and $n$ is even.

\end{enumerate}
\end{conj}

We conclude our paper by
\begin{prob}
Are there any other $\ors$ satisfying Wolstenholme's type theorem
besides those we listed in the paper? More generally is it
possible to find a formula similar to Thm.~\ref{thm:oddlength3}
for arbitrary $\ors$?
\end{prob}

\section*{Appendix: Distribution of Irregular Primes}

Table~\ref{Ta:Irr} in the paper and the following
Table~\ref{Ta:Irr1} give us some evidence to
Conjecture~\ref{conj:dist}.

In Table~\ref{Ta:Irr1} we count the first 30,000 irregular primes
with irregular index $i_p=1$,  first  15,000 irregular primes with
index 2 (producing 30,000 irregular pairs), and all the irregular
primes $<12,000,000$ with index 3 (producing $3\times 9824=29472$
irregular pairs). We denote by $N(k,m)$ the number of irregular
pairs $(p,t)$ satisfying  $p-t\equiv k \pmod 3$ in the top $m$
irregular pairs, $0\le k\le 2$, and
by $P(k,m)$
the percentage of such pairs. We put a subscript $i$ for irregular
primes of index $i$.

\begin{table}[h]
\begin{center}
\begin{tabular}{  ||c|c|c|c|c|c|c|c|c|c|c|| } \hline
$i_p=1,m$&3,000&6,000&9,000&12,000&15,000&18,000&21,000&24,000&27,000&30,000\\
\hline
$N_1(0,m)$&979&1968&2954&4018&5001&5993&6972&7973&8968&9993\\
\hline
$P_1(0,m)$&32.63&32.80&32.82&33.48&33.34&33.29&33.20&33.22&33.22&33.3
\\  \hline
$N_1(1,m)$&1016&2026&3049&4042&5039&6075&7095&8118&9090&10055\\
\hline
$P_1(1,m)$&33.87&33.77&33.88&33.68&33.59&33.75&33.79&33.82&33.67&33.52
\\  \hline
$N_1(2,m)$&1005&2006&2997&3940&4960&5932&6933&7909&8942&9952\\
\hline
$P_1(2,m)$&33.50&33.43&33.30&32.83&33.07&32.96&33.01&32.95&33.12&33.17\\
\hline\hline
$i_p=2,m$&3,000&6,000&9,000&12,000&15,000&18,000&21,000&24,000&27,000&30,000\\
\hline
$N_2(0,m)$&981&2029&3033&4049&5013&6024&7064&8062&9085&10096\\
\hline
$P_2(0,m)$&32.70&33.82&33.70&33.74&33.42&33.47&33.64&33.59&33.65&33.65\\
\hline
$N_2(1,m)$&993&1977&3033&4005&5001&6012&6981&7962&8927&9910\\
\hline
$P_2(1,m)$&33.10&32.95&33.70&33.38&33.34&33.40&33.24&33.18&33.06&33.03\\
\hline
$N_2(2,m)$&1026&1994&2934&3946&4986&5964&6955&7976&8988&9994\\
\hline
$P_2(2,m)$&34.20&33.23&32.60&32.88&33.24&33.13&33.12&33.23&33.29&33.31\\
\hline   \hline
$i_p=3,m$&3,000&6,000&9,000&12,000&15,000&18,000&21,000&24,000&27,000&29,472
\\ \hline
$N_3(0,m)$&997&2027&3065&4102&5124&6091&7087&8067&9046&9889\\
\hline
$P_3(0,m)$&33.23&33.78&34.05&34.18&34.16&33.84&33.74&33.61&33.50&33.55\\
\hline
$N_3(1,m)$&1001&1994&2968&3955&4944&5956&6921&7943&8974&9812\\
\hline
$P_3(1,m)$&33.36&33.23&32.97&32.96&32.96&33.09&32.95&33.09&33.23&33.29\\
\hline
$N_3(2,m)$&1002&1979&2967&3943&4932&5953&6992&7990&8980&9771\\
\hline
$P_3(2,m)$&33.40&32.98&32.96&32.86&32.88&33.07&33.29&33.29&33.26&33.15\\
\hline
\end{tabular}
\caption{Distribution of $p-t\pmod 3$ for irregular pairs $(p,t)$
with $i_p=1,2,3$.} \label{Ta:Irr1}
\end{center}
\end{table}
\begin{table}[h]
\begin{center}
\begin{tabular}{  ||c|c|c|c|c|c|c|c|c|c|| } \hline
$i_p=4,m$ &600&1,200&1,800&2,400&3,000&3,600&4,200&4,800&5,128 \\
\hline $N_4(0,m)$&196&388&569&769&948&1161&1341&1537 & 1645\\
\hline $P_4(0,m)$&32.67&32.33&31.61&32.04&31.60&32.25&31.93&32.02&
32.05\\ \hline $N_4(1,m)$&202&390&615&818&1032&1217&1428&1624 &
1734\\ \hline $P_4(1,m)$&33.67&32.50&34.17&34.08&
34.40&33.81&34.00& 33.83 & 33.79\\ \hline
$N_4(2,m)$&202&422&616&813&1020&1222&1431&1639 & 1749\\ \hline
$P_4(2,m)$&33.67&35.17&34.22&33.88&34.00&33.94&34.07& 34.15 &
34.08\\ \hline
\end{tabular}
\caption{Distribution of $p-t\pmod 3$ for irregular pairs $(p,t)$
with $i_p=4$.} \label{Ta:Irr4}
\end{center}
\end{table}
\begin{table}[h]
\begin{center}
\begin{tabular}{  ||c|c|c|c|c|c|c|c|| } \hline
$i_p=5,m$ &100&200&300&400&500&600&635 \\ \hline
$N_5(0,m)$&31&71&111&145&179&218&230\\ \hline
$P_5(0,m)$&31.00&35.50&37.00&36.25&35.80&36.33&35.94\\ \hline
$N_5(1,m)$&33&58&87&115&154&190&200\\ \hline
$P_5(1,m)$&33.00&29.00&29.00&28.75&30.80&31.67&31.25\\ \hline
$N_5(2,m)$&36&71&102&140&167&192&205\\ \hline
$P_5(2,m)$&36.00&35.50&34.00&35.00&33.40&32.00&32.03\\ \hline
\end{tabular}
\caption{Distribution of $p-t\pmod 3$ for irregular pairs $(p,t)$
with $i_p=5$.} \label{Ta:Irr5}
\end{center}
\end{table}
Between 3 and 12 million there are only 1282 irregular primes with
irregular index 4 (producing 5128 irregular pairs), and 127
irregular primes with irregular index 5 (producing 635 irregular
pairs). We provide the data for them in Table~\ref{Ta:Irr4} and
Table~\ref{Ta:Irr5}. There are 13 irregular primes $<12,000,000$
with irregular index 6, producing 78 pairs. For them we have
$N_6(0,78)=25,$ $N_6(1,78)=27,$ and $N_6(2,78)=26.$ There are
merely 4 irregular primes in the same range with irregular index
7, producing 28 pairs. For these pairs: $N_7(0,28)=10,$
$N_7(1,28)=9,$ and $N_7(2,28)=9,$  which is the best we can hope.

\vskip5cm

\

\newpage

\

\end{document}